\documentclass[a4paper]{article}
\usepackage{amsmath}
    \usepackage{amsxtra}
    \usepackage{amstext}
    \usepackage{amssymb}
    \usepackage{latexsym}

\input xy

\xyoption{all}

\xyoption{poly}

\newcounter{rem}
\newcommand{\rema}{\textbf{Remark 3.1.\arabic{rem}\,}\addtocounter{rem}{1}}

\newtheorem{prop}{Proposition}[section]
\newtheorem{thm}{Theorem}[section]
\newtheorem{thm'}{Theorem}
\newtheorem{cor}{Corollary}[section]
\newtheorem{lem}{Lemma}[section]

\begin{document}
\addtolength{\baselineskip}{+.15\baselineskip} \fontsize{10pt}{13pt}

\begin{center}
\bf \large   GEOMETRICAL DESCRIPTION OF  SMOOTH PROJECTIVE SYMMETRIC
VARIETIES WITH PICARD NUMBER ONE
\end{center}

\

\begin{center}
  \large   Alessandro Ruzzi
\end{center}

\

\begin{abstract}
In \cite{Ru2} we have classified the smooth projective symmetric
$G$-varieties with Picard number one (and $G$ semisimple). In this
work we give a geometrical description of such varieties. In
particular, we determine their group of automorphisms. When this
group, $Aut(X)$, acts non-transitively on $X$, we describe a
$G$-equivariant embedding of the variety $X$ in a homogeneous
variety (with respect to a larger group).
\end{abstract}

{\small\ \ \  \textit{keywords}: Symmetric varieties, Fano
varieties.}

{\small\ \ \ Mathematics Subject Classification 2000: 14M17, 14J45,
14L30}

\

\


A Gorenstein normal algebraic variety $X$ over $\mathbb{C}$ is
called a Fano variety if the anticanonical divisor is ample.  The
Fano surfaces are classically called  Del Pezzo surfaces. The
importance of Fano varieties in the theory of higher dimensional
varieties is similar to the importance of Del Pezzo surfaces in the
theory of surfaces. Moreover  Mori's program predicts that every
uniruled variety is birational to a fiberspace whose general fiber
is a Fano variety (with terminal singularities).

Often it is useful to subdivide the Fano varieties in two kinds: the
Fano varieties with Picard number equal to one and the Fano
varieties whose Picard number is strictly greater of one. For
example, there are many results which give an explicit bound to some
numerical invariants of a Fano variety (depending on the Picard
number and on the dimension of the variety). Often there is an
explicit expression  for the Fano varieties of Picard number equal
to one and another expression for the remaining Fano varieties.

We are mainly interested in the  smooth projective spherical
varieties with Picard number one. The smooth toric (resp.
homogeneous) projective varieties with Picard number one are just
projective spaces (resp. G/P with G simple and P maximal). B.
Pasquier has recently classified the  smooth projective
horospherical varieties with Picard number one (see \cite{P}). In a
previous work we have classified the smooth projective symmetric
$G$-varieties with Picard number one and $G$ semisimple (see
\cite{Ru2}). One can easily show that they are all Fano, because the
canonical bundle cannot be ample. We have also obtained a partial
classification of the smooth Fano complete symmetric varieties with
Picard number strictly greater of one (see \cite{Ru}). Our
classification of the smooth projective symmetric varieties with
Picard number one  is a combinatorial one, so we are naturally
interested to give a geometrical description of such varieties. In
particular, we have proved that, given a symmmetric space $G/H$,
there is at most a smooth completion $X$ of $G/H$ with Picard number
one and $X$ must be projective (see \cite{Ru2}, Theorem 3.1). We
will prove that the automorphism group of a such variety $X$ can act
non-transitively on $X$ only if the rank of $X$ is 2. It would be
interesting to find a reason for such exceptionality of the rank 2
case. Unfortunately, our prove does not explain completely this
fact, because there is a part of the proof that it  is a
case-to-case analysis. The homogeneousness of the rank one varieties
was proved first by Ahiezer in \cite{A}.

More precisely we have proved that:

\begin{thm'}\label{principale not homog}
Let $X$ be a smooth projective completion of a symmetric space $G/H$
with Picard number one (where $G$ is semisimple and simply
connected). Then $Aut(X)$ does not act transitively on $X$ if and
only if: i) the restricted root system has type either $A_{2}$ or
$G_{2}$ and ii) the subgroup $H$ is the subgroup of invariants
$G^{\theta}$.
\end{thm'}

There are six varieties which are not homogenous; their  connected
automorphism group $Aut^{0}(X)$ is isomorphic to $G$ up to isogeny.
These varieties  can be realized as   intersection of hyperplane
sections of a homogeneous projective variety with Picard number one.
Moreover, they are someway related to the exceptional groups; in
particular there are two varieties related to $G_{2}$ and four
varieties related to the magic square of Freudenthal.  The
non-homogeneous varieties with restricted root system of type
$A_{2}$ are obtainable as hyperplane sections of the Legendrian
varieties in the third row of the Freudenthal magic square. Moreover
each one is contained in the the others with bigger dimension. The
completion of $SL_{3}$ was already studied by J. Buczy\'{n}ski (see
\cite{Bu}).

More precisely we will prove the following theorems:

\begin{thm'}\label{principale not homog-explic G2}
Let $G/G^{\theta}$ be a symmetric space whose restricted root system
has type $G_{2}$. We have the following possibilities for the smooth
completion of $G/H$ with Picard number one:

\begin{enumerate}

\item the smooth equivariant completion with Picard number one
of the  symmetric variety $G_{2}/(SL_{2}\times SL_{2})$ of type $G$.
In this case, $Aut(X)=G_{2}$.
The involution $\theta$ can be extended to an involution of
$SO_{7}$,
whose invariant subgroup is $S(O_{3}\times SO_{4})$. The unique
equivariant smooth completion of $SO_{7}/N_{SO_{7}}(S(O_{3}\times
SO_{4}))$ with Picard number one is isomorphic to
$\mathbb{G}_{3}(7)\subset \mathbb{P}^{34}$
and $X$ is the intersection of $\mathbb{G}_{3}(7)$ with a
27-dimensional subspace of $\mathbb{P}^{34}$. 
If we interprete $\mathbb{C}^{7}$ as the subspace of imaginary
elements of the complexified octonions
$\mathbb{O}$, 
then $X$ parametrizes the subspaces $W$ of $\mathbb{C}^{7}$ such
that $W\oplus \mathbb{C}\texttt{1}$ is a subalgebra of $\mathbb{O}$
isomorphic to the complexified quaternions.

\item the smooth equivariant completion with Picard number one
of the  symmetric $(G_{2}\times G_{2})$-variety $G_{2}$. In this
case, $Aut(X)$ is generated by  $ G_{2}\times G_{2}$ and $\theta$.
The involution $\theta$ can be extended to an involution of
$SO_{7}\times SO_{7}$, with invariant subgroup equal to the
diagonal. The unique equivariant smooth completion of $SO_{7} $ with
Picard number one is isomorphic to $\mathbb{IG}_{7}(14)\subset
\mathbb{P}^{63} $ and $X$ is the intersection of
$\mathbb{IG}_{7}(14)$ with a 49-dimensional subspace of
$\mathbb{P}^{63}$.

\end{enumerate}
\end{thm'}

\begin{thm'}\label{principale not homog-explic A2}
Let $G/G^{\theta}$ be a symmetric space whose restricted root system
has type $A_{2}$. We have the following possibilities for the smooth
completion of $G/H$ with Picard number one:
\begin{enumerate}

\item the smooth equivariant completion with Picard number one
of the  symmetric variety $SL_{3}/SO_{3}$ of type $AI$; it is an
hyperplane section of $\mathbb{LG}_{3}(6)$.

\item the smooth equivariant completion with Picard number one
of the  symmetric variety  $SL_{3}$; it is an hyperplane section of
$\mathbb{ G}_{3}(6)$.

\item the smooth equivariant completion with Picard number one
of the  symmetric variety $SL_{6}/Sp_{6}$ of type $AII$; it is an
hyperplane section of $\mathbb{S}_{12}$.

\item the smooth equivariant completion with Picard number one
of the  symmetric variety $E_{6}/F_{4}$ of type $EII$; it is an
hyperplane section of
$E_{7}/P_{7}\equiv\mathbb{G}_{3}(\mathbb{O}^{6})$.
\end{enumerate}

Moreover, $SL_{3}/SO_{3}\subset SL_{3} \subset SL_{6}/Sp_{6}\subset
E_{6}/F_{4}$ and also their smooth completions with Picard number
one are contained nested in each other:

\[ \xymatrix{ \overline{SL_{3}/SO_{3}}\ \ar@{^{(}->}[r]\ar@{^{(}->}[d]& \overline{SL_{3}}
\ \ar@{^{(}->}[r]\ar@{^{(}->}[d]& \overline{SL_{6}/Sp_{6}}\
\ar@{^{(}->}[r]\ar@{^{(}->}[d]&
\overline{E_{6}/F_{4}}\ar@{^{(}->}[d]\\\mathbb{LG}_{3}(6) \
\ar@{^{(}->}[r] & \mathbb{G}_{3}(6)\ \ar@{^{(}->}[r]&
\mathbb{S}_{12}\ \ar@{^{(}->}[r]& \mathbb{G}_{3}(\mathbb{O}^{6})\\
}\] If $G/H$ is different from $SL_{3}$, then the automorphism group
of $X$ is generated by $Aut^{0}(X)$ and $\theta$. If, instead,
$G/H=SL_{3}$ then $Aut^{0}(X)$ has index four in $Aut(X)$.

\end{thm'}

We   give also an explicit description of the  smooth projective
symmetric varieties with Picard number one over which $Aut(X)$ acts
transitively: in particular, we will describe their connected
automorphism group and the immersion of $G/H$ in
$Aut^{0}(X)/Stab_{Aut^{0}(X)}(H/H)$.

In the first section we introduce the notation and recall some
general facts about symmetric varieties. In the second one we cite a
theorem of D. Ahiezer about the homogeneousness of the rank one
varieties. In the third one we study the varieties of rank two which
does not belong to an infinite family; in particular we consider all
the varieties which are not homogeneous. In the forth section we
study the remaining varieties; in particular we will study all the
varieties of rank at least three.

\section{Introduction and notations}\label{intro}

\subsection{Symmetric varieties and colored fans}\label{symm+coloredfan}

In this section we  introduce the necessary  notations. The reader
interested to the embedding theory of spherical varieties can see
\cite{Br3} or \cite{T2}. In \cite{V1} is explained such theory in
the particular case of the symmetric varieties.

\subsubsection{First definitions}\label{1def}

Let $G$ be a connected semisimple algebraic group over $\mathbb{C}$
and let $\theta$ be an involution of $G$. Let $H$ be a closed
subgroup of $G$ such that $G^{\theta}\subset H\subset
N_{G}(G^{\theta})$. We say that   $G/H$ is a symmetric homogeneous
variety. An equivariant embedding of $G/H$ is the data of a
$G$-variety $X$ together with an equivariant open immersion
$G/H\hookrightarrow X$. A normal $G$-variety is called a spherical
variety if it contains a dense orbit under the action of an
arbitrarily chosen Borel subgroup of $G$. One can show that an
equivariant embedding of $G/H$ is a spherical variety if and only if
it is normal (see \cite{dCP1}, Proposition 1.3). In this case we say
that it is a  symmetric variety.   We  say that a subtorus of $G$ is
split if $\theta(t)=t^{-1}$ for all its elements $t$. We say that a
split torus of $G$ of maximal dimension is a maximal split torus and
that a maximal torus containing a maximal split torus is maximally
split. One can prove that any maximally split torus is $\theta$
stable (see \cite{T2}, Lemma 26.5). We fix arbitrarily a maximal
split torus $T^{1}$ and a maximally split torus  $T$ containing
$T^{1}$. Let $R_{G}$ be the root system of $G$ with respect to $T$
and let $R_{G}^{0}$ be the subroot system composed by the roots
fixed by $\theta$. We set $R_{G}^{1}=R_{G}\, \backslash\,
R_{G}^{0}$. We can choose a Borel subgroup $B$ containing $T$ such
that, if $\alpha$ is a positive root in $R_{G}^{1}$, then
$\theta(\alpha)$ is negative (see \cite{dCP1}, Lemma 1.2). One can
prove that $BH$ is dense in $G$ (see \cite{dCP1}, Proposition 1.3).

\subsubsection{Colored fans}\label{colored}

Now, we want to define the colored fan associated to a symmetric
variety. Let $D(G/H)$ be the set of $B$-stable prime divisors of
$G/H$; its elements are called colors. Since $BH/H$ is an affine
open set, the colors are the irreducible components of $(G/H)\
\backslash\, (BH/H)$. We say that a spherical variety is simple if
it contains one closed orbit. Let $X$ be a simple symmetric variety
with closed orbit $Y$. Let $D(X)$ be the subset of $D(G/H)$
consisting of the colors whose closure in $X$ contains $Y$. We say
that $D(X)$ is the set of colors of $X$. To each prime divisor $D$
of $X$, we can associate the normalized discrete valuation $v_{D}$
of $\mathbb{C}(G/H)$ whose ring is the local ring
$\mathcal{O}_{X,D}$. One can prove that $D$ is $G$-stable if and
only if $v_{D}$ is $G$-invariant, i.e. $v_{D}(s\cdot f)=v_{D}(f)$
for each $s\in G$ and $f\in \mathbb{C}(G/H)$. Let $N$ be the  set of
all $G$-invariant valuations of $\mathbb{C}(G/H)$ taking value in
$\mathbb{Z}$ and let $N(X)$ be the set of the valuations associated
to the $G$-stable prime divisors of $X$. Observe that each
irreducible component  of $X\, \backslash\, (G/H)$ has codimension
one, because $G/H$ is affine. Let $S:=T/\,T\cap H\simeq T\cdot
x_{0}$, where $x_{0}=H/H$ denotes the base point of $G/H$. One can
show that the group $\mathbb{C}(G/H)^{(B)}/\,\mathbb{C}^{*}$ is
isomorphic to the character group $\chi(S)$ of $S$ (see \cite{V1},
\S2.3); in particular, it is a free abelian group. We define the
rank of $G/H$ as the rank of $\chi(S)$. We can identify the dual
group
$Hom_{\mathbb{Z}}(\mathbb{C}(G/H)^{(B)}/\,\mathbb{C}^{*},\mathbb{Z})$
with the group  $\chi_{*}(S)$ of one-parameter subgroups of $S$; so
we can identify $\chi_{*}(S)\otimes \mathbb{R}$ with
$Hom_{\mathbb{Z}}(\chi(S),\mathbb{R})$. The  restriction map to
$\mathbb{C}(G/H)^{(B)}/\mathbb{C}^{*}$ is injective over $N$ (see
\cite{Br3}, \S3.1 Corollaire 3), so we can identify $N$ with a
subset of $\chi_{*}(S)\otimes \mathbb{R}$. We say  that $N$ is the
valuation semigroup of $G/H$. Indeed, $N$ is the semigroup
constituted by  the vectors  in the intersection of the lattice
$\chi_{*}(S)$ with an appropriate rational polyhedral convex cone
$\mathcal{C}N$, called  the valuation cone. For each color $D$, we
define $\rho(D)$ as the restriction of $v_{D}$ to $\chi(S)$. In
general, the map $\rho:D(G/H)\rightarrow\chi_{*}(S)\otimes
\mathbb{R}$ is not injective.

Let $C(X)$ be the cone in $\chi_{*}(S)\otimes \mathbb{R}$ generated
by $N(X)$ and $\rho(D(X))$. We denote by $cone(v_{1},...,v_{r})$ the
cone generated by the vectors $v_{1},...,v_{r}$. Given a cone $C$ in
$\chi_{*}(S)\otimes \mathbb{R}$ and a subset $D$ of $D(G/H)$, we say
that  $(C,D)$ is a colored cone if: 1) $C$  is generated by
$\rho(D)$ and by a finite number of vectors in $N$; 2)  the relative
interior of  $C$ intersects $\mathcal{C}N$. The map   $X\rightarrow
(C(X),D(X))$ is a bijection from the set of simple  symmetric
varieties to the set of colored cones (see \cite{Br3}, \S 3.3
Th\'{e}or\`{e}me).

Given a  symmetric variety $\widetilde{X}$ (not necessarily simple),
let $\{Y_{i}\}_{i\in I}$ be \vspace{0.1 mm} the set of $G$-orbits.
Observe that $\widetilde{X}$ contains a finite number of $G$-orbits,
thus $\widetilde{X}_{i}:=\{x\in \widetilde{X}\ |\  \overline{G\cdot
x}\supset Y_{i}\}$ is open in $\widetilde{X}$ and is a simple
symmetric variety whose closed orbit is $Y_{i}$. We define
$D(\widetilde{X})$ as the set $\bigcup_{i\in
I}D(\widetilde{X}_{i})$. The family
$\{(C(\widetilde{X}_{i}),D(\widetilde{X}_{i}))\}_{i\in I}$ is
\vspace{0.2 mm} called the colored fan of $\widetilde{X}$ and
determines completely $\widetilde{X}$ (see \cite{Br3}, \S 3.4
Th\'{e}or\`{e}me 1). Moreover $\widetilde{X}$ is complete if and
only if $\mathcal{C}N\subset\bigcup_{i\in I}C(\widetilde{X}_{i}) $
(see \cite{Br3}, \S 3.4 Th\'{e}orem\`{e} 2).

\subsubsection{Restricted root system}\label{restrictedrootsystem}

To describe  the sets $N$ and $\rho(D(G/H))$, we need to associate a
root system  to $G/H$. The subgroup $\chi(S)$ of $\chi(T^{1})$ has
finite index, so we can identify $\chi(T^{1} )\otimes\,\mathbb{R}$
with $\chi(S)\otimes\,\mathbb{R}$. Because $T$ is $\theta$-stable,
$\theta$ induces an involution of $\chi(T)\otimes\,\mathbb{R}$ which
we call again $\theta$. The inclusion $T^{1}\subset T$ induces an
isomorphism of $\chi(T^{1})\otimes\,\mathbb{R}$ with the
$(-1)$-eigenspace of $\chi(T)\otimes\,\mathbb{R}$ under the action
of $\theta$ (see \cite{T2}, \S 26). Denote  by $W_{G}$ the Weyl
group of $G$ (with respect to $T$) and let $(\, \cdot ,\cdot)$ be
the Killing form over $span_{\mathbb{R}}(R_{G})$. We denote with the
same symbol the restriction of  $(\, \cdot ,\cdot)$ to
$\chi(T^{1})\otimes\,\mathbb{R}$, thus we can identify
$\chi(T^{1})\otimes\,\mathbb{R}$ with its dual
$\chi_{*}(T^{1})\otimes\,\mathbb{R}$.  The set $R_{G,\theta}:=
\{\beta-\theta(\beta)\ |\ \beta\in R_{G}\}\backslash\,\{0\}$ is a
root system in $\chi(S)\otimes\mathbb{R}$ (see \cite{V1}, \S 2.3
Lemme), which we call the restricted root system of $(G,\theta)$; we
call the non zero $\beta-\theta(\beta)$ the restricted roots.
Usually we denote by $\beta$ (respectively by $\alpha$)  a root of
$R_{G}$ (respectively of $R_{G,\theta}$); often we denote by
$\varpi$ (respectively by $\omega$)  a weight of $R_{G}$
(respectively of $R_{G,\theta}$). In particular, we denote by
$\varpi_{1},...,\varpi_{n}$ the fundamental weight of $R_{G}$ (we
have chosen the basis of $R_{G}$ associated to $B$). Notice however
that the weights of $R_{G,\theta}$ are weights of $R_{G}$. The
involution $\iota:=-\varpi_{0}\cdot\theta$ of $\chi(T)$ preserves
the set of simple roots; moreover $\iota$ coincides with $-\theta$
modulo the lattice generated by $R_{G}^{0}$ (see \cite{T2}, p. 169).
Here $\varpi_{0}$ is the longest element of the Weyl group of
$R_{G}^{0}$.  We denote by $\alpha_{1},...,\alpha_{s}$ the elements
of the basis $\{\beta-\theta(\beta)\, |\, \beta\in R_{G}$
simple$\}\,\backslash\, \{0\}$   of $R_{G,\theta}$. Let $b_{i}$ be
equal to $\frac{1}{2}$ if $2\alpha_{i}$ belongs to $R_{G,\theta}$
and equal to one otherwise; for each $i$ we define
$\alpha_{i}^{\vee}$ as the coroot
$\frac{2b_{i}}{(\alpha_{i},\alpha_{i})}\alpha_{i}$. The set
$\{\alpha^{\vee}_{1} ,...,\alpha^{\vee}_{s}\}$ is a basis of the
dual root system $R^{\vee}_{G,\theta}$, namely the root system
composed by the coroots of the restricted roots. We call the
elements of $R^{\vee}_{G,\theta}$ the restricted coroots.  Let
$\omega_{1},...,\omega_{s}$  be the fundamental weights of $R
_{G,\theta}$  with respect to $\{\alpha_{1} ,...,\alpha_{s}\}$ and
let $\omega^{\vee}_{1},...,\omega^{\vee}_{s}$  be the fundamental
weights of $R^{\vee}_{G,\theta}$ with respect to
$\{\alpha^{\vee}_{1} ,...,\alpha^{\vee}_{s}\}$. Let $C^{+}$ be the
positive closed Weyl chamber of $\chi(S)\otimes\mathbb{R}$; we call
$-C^{+}$ the negative Weyl chamber. Given a dominant weight
$\lambda$ of $G$, we denote by $V(\lambda)$ the irreducible
representation of highest weight $\lambda$.

We want to give a description of the  weight lattice of
$R_{G,\theta}$. We say that a dominant weight $\varpi\in \chi(T)$ is
a spherical weight if $V(\varpi)$ contains a non-zero vector fixed
by $G^{\theta}$. In this case, $V^{G^{\theta}}$ is one-dimensional
and $\theta(\varpi)=-\varpi$. Thus, we can think $\varpi$ as a
vector in $\chi(S)\otimes\mathbb{R}$. One can show that  the lattice
generated by the spherical weights coincides with the weight lattice
of $R_{G,\theta}$. See \cite{CM}, Theorem 2.3 or \cite{T2},
Proposition 26.4 for an explicit description of the spherical
weights. Such description implies that the set of dominant weights
of $R_{G,\theta}$ is the set of spherical weights and that $C^{+}$
is the intersection of $\chi(S)\otimes\mathbb{R}$ with the positive
closed Weyl chamber of the root system $R_{G}$.

\subsubsection{The sets $N$ and $D(G/H)$}\label{n+dg/h}

The set $N$ is equal to $-C^{+}\cap\,\chi_{*}(S)$; in particular, it
consists of the lattice vectors of the rational polyhedral convex
cone $\mathcal{C}N=-C^{+}$. The set $\rho(D(G/H))$ is equal to
$\{\alpha^{\vee}_{1} ,...,\alpha^{\vee}_{s}\}$ (see \cite{V1}, \S
2.4, Proposition 1 and Proposition 2). For each $i$, the fibre
$\rho^{-1}(\alpha_{i}^{\vee})$ contains  at most 2 colors (see
\cite{V1}, \S2.4, Proposition 1). We say that a simple restricted
root $\alpha_{i_{0}}$ is exceptional if there are two distinct
simple roots $\beta_{i_{1}}$ and $\beta_{i_{2}}$ in $R_{G}$ such
that: 1) $\beta_{i_{1}}-\theta(\beta_{i_{1}})=
\beta_{i_{2}}-\theta(\beta_{i_{2}}) = \alpha_{i_{0}}$; 2) either
$\theta(\beta_{i_{1}})\neq -\beta_{i_{2}}$ or
$\theta(\beta_{i_{1}})= -\beta_{i_{2}}$ and $(\beta_{i_{1}},
\beta_{i_{2}})\neq0$. In this case we say that also
$\alpha_{i_{0}}^{\vee}$, $\theta$ and all the equivariant embeddings
of $G/H$ are exceptional. If $G/H$ is exceptional, then $\rho$ is
not injective.
We say that $G/H$ contains a Hermitian factor if the center of
$[G,G]^{\theta}$ has positive dimension.  If  $G/H$ does not contain
a Hermitian factor, then $\rho$ is injective (see \cite{V1}, \S 2.4,
Proposition 1). If $\rho$ is injective, we denote by
$D_{\alpha^{\vee}}$ the unique color contained in
$\rho^{-1}(\alpha^{\vee})$.

\subsubsection{Local description of  symmetric
varieties}\label{localdescription}

Let $X$ be a  symmetric variety 
and let $Y$ be a closed orbit of $X$. One can show that there is a
unique  $B$-stable affine open set $X_{B}$ that intersects $Y$ and
is minimal for this property. Moreover, the complement
$X\,\backslash\, X_{B}$ is the union of  the $B$-stable prime
divisors   not containing $Y$ (see \cite{Br3}, \S 2.2, Proposition).
The stabilizer $P$ of $X_{B}$ is a parabolic subgroup containing
$B$. Suppose $X$ smooth or non-exceptional and let $L$ be the
standard Levi subgroup of $P$. Then  $L$ is $\theta$-stable,
$\theta(P)=P^{\,\text{-}}$ and the $P$-variety $X_{B}$ is the
product $R_{u}P\times Z$ of the unipotent radical of $P$ and of an
affine $L$-symmetric variety $Z$ containing $x_{0}$  (see
\cite{Br3}, \S 2.3, Th\'{e}or\`{e}me and \cite{Ru2}, Lemma 2.1). If
$Y$ is projective, then $Z$ contains a $L$-fixed point whose
stabilizer in $G$ is $P^{\,\text{-}}$. Given a root $\beta$, let
$U_{\beta}$ be the unipotent one-dimensional subgroup of $G$
corresponding to $\beta$. Given $\mu\in
\chi_{*}(T)\otimes\mathbb{Q}\equiv \chi (T)\otimes\mathbb{Q}$, we
denote by $P(\mu)$ the parabolic subgroup of $G$ generated by $T$
and by the subgroups $U_{\beta}$ corresponding to the roots $\beta$
such that $(\beta,\mu)\ \geq 0$. Given a parabolic subgroup
$P=P(\mu)$, sometimes we denote by $P^{\,\text{-}}$ the opposite
standard  parabolic subgroup, namely $P(-\mu)$.

\subsection{Colored fans of smooth complete symmetric varieties with Picard number
one}\label{smoothcoloredfan}

In \cite{Ru2} we have proved the following combinatorial
classification of the smooth complete symmetric varieties with
Picard number one.  \emph{Unless otherwise specified, we always
suppose $G$ simply connected.} Notice that this assumption  is not
restrictive   (see, for example, \cite{V1}, \S2.1).

\begin{thm}\label{classif} Let $G$ be a semisimple, simply connected group and let
$G/H$ be a homogeneous symmetric variety.  Suppose that there is a
smooth, complete embedding $X$ of $G/H$ with Picard number one.
Then:
\begin{itemize}
\item Fixed $G/H$, there is, up to  equivariant isomorphism, at most one  embedding with
the previous properties. Moreover,  it is  projective and contains
at most two closed orbits.

\item The number of colors of $G/H$ is equal to the rank $l$ of $G/H$; in
particular there are no exceptional roots.
\item We have the following classification
depending on the type of the restricted root system $R_{G,\theta}$ :
\begin{enumerate}
\item[(i)] If $R_{G,\theta}$  has type $A_{1}\times A_{1}$, then $\chi(S)$ has basis
$\{2\omega_{1},\omega_{1}+\omega_{2}\}$; in particular $H$ has index
two in $N_{G}(G^{\theta})$. Moreover, $X$ has two closed orbits; the
maximal colored cones of the colored fan of $X$ are
$(cone(\alpha^{\vee}_{1},-\omega^{\vee}_{1}-\omega^{\vee}_{2}),$ $
\{D_{\alpha^{\vee}_{1}}\})$ and
$(cone(\alpha^{\vee}_{2},-\omega^{\vee}_{1}-\omega^{\vee}_{2}),$ $
\{D_{\alpha^{\vee}_{2}}\})$.

\item[(ii)] If $l=1$, then $G/H$ can   be isomorphic neither to
$SL_{n+1}/S(L_{1}\times L_{n})$, nor to $SL_{2}/SO_{2}$. With such
hypothesis, $G/H$ has a  unique non trivial embedding which is
simple, projective, smooth and with Picard number one.
\item[(iii)] If $R_{G,\theta}$  has type $A_{l}$ with $l>1$, we have the following
possibilities:
\begin{itemize}
\item[(a)] $H=N_{G}(G^{\theta})$ and $X$ is simple. In this case $X$ is
associated either to the colored cone
$(cone(\alpha^{\vee}_{1},...,\alpha^{\vee}_{l-1},-\omega^{\vee}_{1}),
\{D_{\alpha^{\vee}_{1}},...,D_{\alpha^{\vee}_{l-1}}\})$ or to the
one
$(cone(\alpha^{\vee}_{2},...,\alpha^{\vee}_{l},-\omega^{\vee}_{l}),
\{D_{\alpha^{\vee}_{2}},...,D_{\alpha^{\vee}_{l }}\})$;
\item[(b)] $H=G^{\theta}$ and $l=2$. In this case $X$ has two closed orbits. The
maximal colored cones of the colored fan of $X$ are
$(cone(\alpha^{\vee}_{1},-\omega^{\vee}_{1}-\omega^{\vee}_{2}),
\{D_{\alpha^{\vee}_{1}} \})$ and
$(cone(\alpha^{\vee}_{2},-\omega^{\vee}_{1}-\omega^{\vee}_{2}),
\{D_{\alpha^{\vee}_{2}}\})$.
\end{itemize}
\item[(iv)] If $R_{G,\theta}$  has type $B_{2}$, then $X$ is simple and  we have the following
possibilities:

\begin{itemize}
\item[(a)] $H=N_{G}(G^{\theta})$ and  $X$ is
associated   to  $(cone(\alpha^{\vee}_{1},-\omega^{\vee}_{1}),$ $
\{D_{\alpha^{\vee}_{1}}\})$;
\item[(b)] $H=G^{\theta}$ and $X$ is
associated   to  $(cone(\alpha^{\vee}_{2},-\omega^{\vee}_{2}),$
$\{D_{\alpha^{\vee}_{2}}\})$. Moreover  $G/H$ cannot be Hermitian.
\end{itemize}

\item[(v)] If $R_{G,\theta}$  has type $B_{l}$ with $l>2$, then $H=N_{G}(G^{\theta})$,  $X$ is
simple and  is associated  to
$(cone(\alpha^{\vee}_{1},...,\alpha^{\vee}_{l-1},-\omega^{\vee}_{1}),
\{D_{\alpha^{\vee}_{1}},...,$ $D_{\alpha^{\vee}_{l-1}}\})$.

\item[(vi)] If $R_{G,\theta}$  has type $C_{l}$,   then $H= G^{\theta} $, $X$ is simple and
corresponds to
$(cone(\alpha^{\vee}_{1},...,\alpha^{\vee}_{l-1},-\omega^{\vee}_{1}),
\{D_{\alpha^{\vee}_{1}},...,D_{\alpha^{\vee}_{l-1}}\})$. Moreover,
$G/H$ cannot be Hermitian.

\item[(vii)] If $R_{G,\theta}$  has type $BC_{l}$ with $l>1$, then $H=N_{G}(G^{\theta})= G^{\theta}$,  $X$ is
simple and corresponds  to
$(cone(\alpha^{\vee}_{1},...,\alpha^{\vee}_{l-1},-\omega^{\vee}_{1}),$
$ \{D_{\alpha^{\vee}_{1}},...,$ $D_{\alpha^{\vee}_{l-1}}\})$.
\item[(viii)] If $R_{G,\theta}$  has type $D_{l}$ with $l>4$, then $\chi_{*}(S)$ is freely generated by
$\omega^{\vee}_{1},...,\omega^{\vee}_{l-2},\omega^{\vee}_{l-1}+\omega^{\vee}_{l},2\omega^{\vee}_{l}$;
in particular $H$ has index two in $N_{G}(G^{\theta})$. $X$ has two
closed orbits; the maximal colored cones of the colored fan of $X$
are
$(cone(\alpha^{\vee}_{1},...,\alpha^{\vee}_{l-1},-\omega^{\vee}_{1}),$
$ \{D_{\alpha^{\vee}_{1}},...,D_{\alpha^{\vee}_{l-1}}\})$ and
$(cone(\alpha^{\vee}_{1},...,$ $ \alpha^{\vee}_{l-2},$
$\alpha^{\vee}_{l}, -\omega^{\vee}_{1}),$
$\{D_{\alpha^{\vee}_{1}},...,$ $D_{\alpha^{\vee}_{l-2}},$
$D_{\alpha^{\vee}_{l }}\})$.
\item[(ix)] If $R_{G,\theta}$  has type $D_{4}$, then $H$ has index two in $N_{G}(G^{\theta})$ and $X$
has two closed orbits. If
$\chi_{*}(S)=\mathbb{Z}\omega^{\vee}_{i}\oplus
\mathbb{Z}\omega^{\vee}_{2}\oplus
\mathbb{Z}(\omega^{\vee}_{j}+\omega^{\vee}_{k})\oplus
\mathbb{Z}2\omega^{\vee}_{k}$, then the maximal colored cones of the
colored fan of  $X$ are
$(cone(\alpha^{\vee}_{i},\alpha^{\vee}_{2},\alpha^{\vee}_{j},$
$-\omega^{\vee}_{i}),$ $
\{D_{\alpha^{\vee}_{i}},D_{\alpha^{\vee}_{2}},D_{\alpha^{\vee}_{j}}\})$
and
$(cone(\alpha^{\vee}_{i},\alpha^{\vee}_{2},\alpha^{\vee}_{k},-\omega^{\vee}_{i}),$
$ \{D_{\alpha^{\vee}_{i}},D_{\alpha^{\vee}_{2}},$
$D_{\alpha^{\vee}_{k}}\})$.
\item[(x)] If $R_{G,\theta}$  has type $G_{2}$ then $H=N_{G}(G^{\theta})= G^{\theta}$,  $X$ is
simple and  is associated  to
$(cone(\alpha^{\vee}_{2},-\omega^{\vee}_{2}),
\{D_{\alpha^{\vee}_{2}}\})$.
\item[(xi)]
If the type of $R_{G,\theta}$ is different from the previous ones,
then there is no a variety $X$ with the requested properties.

\end{enumerate}
\end{itemize}
\end{thm}

\subsection{Description of some exceptional groups via composition
algebras}\label{exceptionalgroup}

We will need a description of some exceptional groups via complex
composition algebras and Jordan algebras.  The interested reader can
see \cite{LM} and \cite{Ad} for a  detailed exposition of the facts
which we recall here.

Let $\mathbb{A}$ be a complex composition algebra, i.e.
$\mathbb{A}=\mathbb{A}_{\mathbb{R}}\otimes_{\mathbb{R}}\mathbb{C}$
where $\mathbb{A}_{\mathbb{R}}$ is a real division algebra (namely
$\mathbb{A}_{\mathbb{R}}=\mathbb{R}$, $\mathbb{C}$, $\mathbb{H}$ or
$\mathbb{O}$). If $a\in \mathbb{A}$, let $\bar{a}$ be its conjugate;
we denote by $Im \mathbb{A}$ the subspace of  pure imaginary
element, i.e. the elements $a$ such that $\bar{a}=-a$. Let
$\mathcal{J}_{3}(\mathbb{A})$ be the space of $\mathbb{A}$-Hermitian
matrices of order three, with coefficients in $\mathbb{A}$:

\[\mathcal{J}_{3}(\mathbb{A})=\left\{\left(\begin{matrix}
r_{1}&\bar{x}_{3}&\bar{x}_{2}\\
x_{3}&r_{2}      &\bar{x}_{1}\\
x_{2}&x_{1}      &r_{3}
\end{matrix}\right), r_{i}\in\mathbb{C},x_{i}\in\mathbb{A}\right\}.\]

$\mathcal{J}_{3}(\mathbb{A})$  has the structure of a Jordan algebra
with the multiplication  $A\circ B:=\frac{1}{2}(AB+BA)$, where $AB$
is the usual matrix multiplication. There is a well defined cubic
form on $\mathcal{J}_{3}(\mathbb{A})$, which we call the
determinant. Given $P\in \mathcal{J}_{3}(\mathbb{A})$, its comatrix
is defined by
\[com\,P=P^{2}-(trace\,P)P+\frac{1}{2}((trace\,P)^{2}-trace\,P^{2})I \]
and characterized by the identity $com(P)P=det(P)I$. In particular,
$P$ is invertible if and only if $det(P)$ is different from 0.

Let $SL_{3}(\mathbb{A})\subset
GL_{\mathbb{C}}(\mathcal{J}_{3}(\mathbb{A}))$ be subgroup preserving
the determinant; $\mathcal{J}_{3}(\mathbb{A})$ is an irreducible
$SL_{3}(\mathbb{A})$ representation. We let $SO_{3}(\mathbb{A})$
denote the group of complex linear transformations preserving the
Jordan multiplication; it is also the subgroup of
$SL_{3}(\mathbb{A})$ preserving the quadratic form
$Q(A)=trace(A^{2})$.

We call $\mathcal{Z}_{2}(\mathbb{A}):=\mathbb{C}\oplus
\mathcal{J}_{3}(\mathbb{A})\oplus\mathcal{J}_{3}(\mathbb{A})^{*}\oplus\mathbb{C}^{*}$
the space of Zorn matrices. The space
$\mathfrak{s}\mathfrak{p}_{6}(\mathbb{A}):=\mathbb{C}^{*}\oplus
\mathcal{J}_{3}(\mathbb{A})^{*}\oplus(Lie(SL_{3}(\mathbb{A}))+\mathbb{C})
\oplus\mathcal{J}_{3}(\mathbb{A})\oplus\mathbb{C}$ has a structure
of Lie algebra and $\mathcal{Z}_{2}$ has a natural structure of
(simple) $\mathfrak{s}\mathfrak{p}_{6}(\mathbb{A})$-module. There is
a unique closed connected subgroup of
$GL_{\mathbb{A}}(\mathcal{Z}_{2}(\mathbb{A}))$ with Lie algebra
$\mathfrak{s}\mathfrak{p}_{6}(\mathbb{A})$; we denote it by
$Sp_{6}(\mathbb{A})$. Moreover, there is a
$Sp_{6}(\mathbb{A})$-invariant symplectic form on
$\mathcal{Z}_{2}(\mathbb{A})$.

The closed $Sp_{6}(\mathbb{A})$-orbit in
$\mathbb{P}(\mathcal{Z}_{2}(\mathbb{A}))$  is the image of the
$Sp_{6}(\mathbb{A})$-equivariant rational map:

\[\begin{matrix}
\phi:\mathbb{P}(\mathbb{C}\oplus\mathcal{J}_{3}(\mathbb{A}))&\dashrightarrow&
 \mathbb{P}(\mathcal{Z}_{2}(\mathbb{A}))\\
 (x,P)&\rightarrow &(x^{3},x^{2}P,x\,com(P),det(P)).\end{matrix}\]
Furthermore, if $\mathbb{C}$ is interpreted as a space of diagonal
matrix (in $\mathcal{J}_{3}(\mathbb{A})$) and $(I,P)$ is interpreted
as a matrix of three row vectors in $\mathbb{A}^{6}$, then the
previous map is the usual Plucker map. The condition $P\in
\mathcal{J}_{3}(\mathbb{A})$ can be interpreted as the fact that the
three vectors defined by the matrix $(I,P)$ are orthogonal with
respect to the Hermitian symplectic two-form $w(x,y)=^{t}xJ\bar{y}$,
where $J=\left(\begin{smallmatrix}0&I\\
\text{-}  I&0\end{smallmatrix}\right)$.
Therefore, it is natural to see the closed
$Sp_{6}(\mathbb{A})$-orbit in
$\mathbb{P}(\mathcal{Z}_{2}(\mathbb{A}))$ as the Grassmannian
$\mathbb{LG}(\mathbb{A}^{3},\mathbb{A}^{6})$ of symplectic 3-planes
in $\mathbb{A}^{6}$.

Explicitly, we have the following possibilities: 1) if
$\mathbb{A}_{R}=\mathbb{R}$ then $SL_{3}(\mathbb{A})$ is $SL_{3}$,
$SO_{3}(\mathbb{A})$ is $SO_{3}$ and $Sp_{6}(\mathbb{A})$ is
$Sp_{6}$; 2) if $\mathbb{A}_{R}=\mathbb{C}$ then
$SL_{3}(\mathbb{A})$ is $SL_{3}\times SL_{3}$, $SO_{3}(\mathbb{A})$
is $SL_{3}$ and $Sp_{6}(\mathbb{A})$ is $SL_{6}$; 3) if
$\mathbb{A}_{R}=\mathbb{H}$ then $SL_{3}(\mathbb{A})$ is $SL_{6}$,
$SO_{3}(\mathbb{A})$ is $Sp_{6}$ and $Sp_{6}(\mathbb{A})$ is
$Spin_{12}$; 4) if $\mathbb{A}_{R}=\mathbb{O}$ then
$SL_{3}(\mathbb{A})$ is $E_{6}$, $SO_{3}(\mathbb{A})$ is $F_{4}$ and
$Sp_{6}(\mathbb{A})$ is $E_{7}$.

\subsection{Generalized flag varieties}\label{flag}

We conclude this section with a  description of the projective
homogeneous varieties with Picard number one for the classic groups.
Let $V$ be a $n$-dimensional vector space, we will denote by
$\mathbb{G}_{m}(V)$ the \textit{Grassmannian of the $m$-dimensional
subspace} of $V$. Let $q$ be a non-degenerate symmetric bilinear
form on $V$ and let $SO(V,q)$ be the corresponding special
orthogonal group. We will say that a subspace of $V$ is isotropic if
the restriction of $q$ to it is zero. Fix an integer $m$ such that
$2m\leq r$ and let $\mathbb{I}\mathbb{G}_{m}(V)\subset
\mathbb{G}_{m}(V)$ be the algebraic subvariety whose points are
identified with isotropic $m$-dimensional subspaces of $V$. The
group $SO(V,q)$ acts on $\mathbb{I}\mathbb{G}_{m}(V)$ in the natural
manner and the action is transitive if $2m< n$. Instead, if $2m= n$,
$\mathbb{I}\mathbb{G}_{m}(V)$ consists of two isomorphic
$SO(V,q)$-orbits (each of them being a connected component of
$\mathbb{I}\mathbb{G}_{m}(V)$); we  denote a such orbit by
$\mathbb{S}_{m}(V)$ (or by $\mathbb{S}_{m} $ if $V=\mathbb{C}^{2m}$
and $q$ is the standard symmetric bilinear form). The variety
$\mathbb{I}\mathbb{G}_{m}(V)$ is called the \textit{isotropic
Grassmannian of $m$-dimensional isotropic subspace}, while
$\mathbb{S}_{m}(V)$ is called the \textit{spinorial variety of order
$m$}. In an analogous manner, let $q'$ be a non-degenerate
skew-symmetric bilinear form on $V$ and fix an integer $m$ such that
$2m\leq n$. We denote by $\mathbb{L}\mathbb{G}_{m}(V)$ the algebraic
subvariety of $\mathbb{G}_{m}(V)$ whose points are identified with
isotropic $m$-dimensional subspaces of $V$;
$\mathbb{L}\mathbb{G}_{m}(V)$ is called the \textit{Lagrangian
Grassmannian}.

\section{Varieties of rank one}\label{sect-rank1}

In this section we  describe the varieties  with rank one. For any
homogenous symmetric variety $G/H$ with rank one there is a unique
equivariant completion  $X$. It is smooth, projective and has Picard
number at most two. Moreover, $X$  has  exactly two orbits: an open
orbit and a closed orbit of codimension one.  Thus $X$ is a
homogeneous variety by the following theorem due to D. Ahiezer.

\begin{thm}[see \cite{A}, Theorem 4]\label{rank 1} Let $G$ be a semisimple group and
let $H$ be a closed reductive subgroup. Let $X$ be an equivariant
(normal) completion of $G/H$ such that $X\,\setminus\, (G/H)$ is a
$G$-orbit (of codimension one).  Then $X$ is a homogeneous space for
a larger group. If $G/H$ is a symmetric variety and $G$ acts
effectively on $G/H$, we have the following possibilities:
\begin{enumerate}
\item $G$ is $SL_{2}\times SL_{2}$, $H$ is $SL_{2}$ and $X$ is 
$\{(x,t)\,|\, det\,(
x)=t^{2}\}\subset \mathbb{P}(S^{2}\mathbb{C}^{2}\oplus \mathbb{C})$;

\item $G$ is $PSL_{2}\times PSL_{2}$, $H$ is $PSL_{2}$  and $X$ is 
$\mathbb{P}(S^{2}\mathbb{C}^{2})$;

\item $G$ is $SL_{n}$, $H$ is $GL_{n-1}$, $\theta$ has type $AIV$ (or $AI$ if $n=2$) and $X$ is 
$\mathbb{P}(\mathbb{C}^{n})\times \mathbb{P}((\mathbb{C}^{n})^{*})$;

\item $G$ is $PSL_{2}$, $H$ is $PSO_{2}$, $\theta$ has type $AI$  and $X$ is 
$\mathbb{P}(\mathfrak{sl}_{2})$;

\item $G$ is $Sp_{2n}$, $H$ is $Sp_{2}\times Sp_{2n-2}$, $\theta$ has type $CII$
and $X$ is 
$\mathbb{G}_{2}(2n)$;

\item $G$ is $SO_{n}$,  $H$ is $SO_{n-1}$, $\theta$ has type $BII$ or $DII$
and $X$ is 
$\{(x,t)\,|\, q(x,x)=t^{2}\}\subset \mathbb{P}(\mathbb{C}^{n}\oplus
\mathbb{C})$, where $q$ is the standard symmetric bilinear form;

\item $G$ is $SO_{n}$, $H$ is $S(O_{1}\times O_{n-1})$, $\theta$ has type $BII$ or $DII$,
and $X$ is 
$ \mathbb{P}(\mathbb{C}^{n})$;

\item $G$ is $F_{4}$,  $H$ is $Spin_{9}$, $\theta$ has type $FII$ and $X$ is 
$E_{6}/P_{\omega_{1}}\equiv \mathbb{P}^{2}\mathbb{O}$;

\end{enumerate}
\end{thm}

\section{Varieties of rank two}\label{rank2}

\textit{In the following of this work we always suppose that $G/H$
has rank strictly greater than one.} In this section we describe the
varieties of rank  two which do not belong to an infinite family.
Explicitly, we consider completion of the following homogenous
varieties: 1) the symmetric variety $SL_{4}/N_{SL_{4}}(SL_{2}\times
SL_{2})$ of type $AIII$; 2) the symmetric  $(SL_{3}\times
SL_{3})$-variety $SL_{3}$; 3) the symmetric variety $SL_{3}/SO_{3}$
of type $AI$; 4) the symmetric variety $SL_{4}/Sp_{4}$ of type
$AII$;  5) the symmetric variety $E_{6}/F_{4}$ of type $EIV$ (with
$E_{6}$ simply connected); 6) the
 symmetric variety $E_{6}/N_{E_{6}}(F_{4})$ of type $EIV$; 7) the
symmetric variety $Sp_{4}/(Sp_{2}\times Sp_{2})$ of type $CII$; 8)
the symmetric variety $G_{2}$; 9) the symmetric variety
$G_{2}/(Sl_{2}\times Sl_{2})$ of type $G$; 10) the symmetric
varieties whose  restricted root system  has type $A_{1}\times
A_{1}$. In the last case $G/N_{G}(G^{\theta})$ is isomorphic to
$SO_{n}/N_{SO_{n}}(SO_{n-1})\times SO_{m}/N_{SO_{m}}(SO_{m-1})$;
$n,m$ are strictly greater than two  and $H$ has index two in
$N_{G}(G^{\theta})$.

In this section we do not consider the completions of the symmetric
varieties $Spin_{5}$ and $SO_{8}/N_{SO_{8}}(GL_{4})$  because: 1)
$Spin_{5}$ is isomorphic to $Sp_{4}$ and we will study it together
with the completion of $Sp_{2n}$; 2) $SO_{8}/$ $N_{SO_{8}}(GL_{4})$
is isomorphic to $SO_{8}/N_{SO_{8}}(SO_{2}\times SO_{6 })$ and we
will study it together with the completion of
$SO_{n}/N_{SO_{n}}(SO_{2}\times SO_{n-2})$.

\subsection{Some general considerations}\label{general}

We begin with some general considerations; in particular we do not
suppose yet that $X$ has rank two. We  prove that a   smooth
projective symmetric variety $X$ with Picard number one is a
homogeneous variety if and only if
$H^{0}(G/P^{\,\text{-}},N_{G/P^{\,\text{-}},\,X})\neq0$ for an
appropriate closed orbit $G/P^{\,\text{-}}$  of $X$. Furthermore, we
explain how to study the automorphism group of $X$ (if such variety
is not homogenous).

\begin{lem}\label{smooth subvarieties} Let $X$ be a smooth projective  symmetric variety
with Picard number one and let $Y$ be a proper $G$-stable  closed
subvariety of $X$. If $Y$ is smooth, then it is a $G$-orbit.
\end{lem}

{\em Proof.} Recall that the minimal $B$-stable affine open set
$X_{B}$ which intersects $Y$ is a product $P_{u}\times Z$, where $Z$
is an affine $L$-symmetric variety. One can easily show that $Z$ is
an irreducible $L$-representation (see \cite{Ru2}, pages 9 and 17 or
\cite{Ru2}, Theorem 2.2). Observe that $Y$ is smooth if and only if
$Y\cap Z$ is smooth. But $Y\cap Z$ is a cone because the center of
$L$ acts non-trivially on it. Thus $Y\cap Z$ is smooth if and only
if it is a subrepresentation of $Z$. Since $Z$ is irreducible,
$Y\cap Z$ is smooth if and only if it is a point. In such case $Y$
is a closed orbit.
$\square$

We  use the previous lemma only to prove following corollary.

\begin{cor}\label{closed orbit} Let $X$ be a smooth projective symmetric
variety with Picard number one. If  $Aut^{0}(X)$ does not act
transitively over $X$, then $Aut^{0}(X)$ stabilizes a closed
$G$-orbit.
\end{cor}

{\em Proof.} Let $Y$ be a minimal closed $G$-stable subvariety
stabilized by $Aut^{0}(X)$. Suppose by contradiction that $Y$ is not
a $G$-orbit, then it is singular by the Lemma~\ref{smooth
subvarieties}. Thus the singular locus of $Y$ is not empty and is
stabilized by $Aut^{0}(X)$. In particular,  its irreducible
components   are  $Aut^{0}(X)$-stable, closed subvarieties of $X$,
properly contained in $Y$; a contradiction. $\square$

A closed 
orbit $G/P^{\,\text{-}}$ of $X$ is stabilized by
$Aut^{0}(X)$ if and only if $H^{0}(G/P^{\,\text{-}},$ $
N_{G/P^{\,\text{-}},X})=0$, where $N_{G/P^{\,\text{-}},X}$ is the
normal bundle to $G/P^{\,\text{-}}$ (see \cite{A2}, $\S2.3$). Recall
that $X_{B}$ is $P_{u}\times Z$; m%
oreover, the intersection of $G/P^{\,\text{-}}$ with $X_{B}$ is
$P_{u}\times \{0\}$. Thus we can identify $Z$ with the fibre of the
normal bundle $N_{G/P^{\,\text{-}},\,X}$ over
$P^{\,\text{-}}/P^{\,\text{-}}$. We know by the Borel-Weil theorem
(see \cite{A2}, \S4.3) that $H^{0}(G/P^{\,\text{-}},$
$N_{G/P^{\,\text{-}},\,X})=0$ if and only if the highest weight of
$Z$ is not dominant.

Now, we want to explain how to calculate the highest weight $\omega$
of $Z$ when the  rank of $G/H$ is two. The center of $L$ has
dimension one and $R_{[L,L],\theta}$ has rank one; moreover, if
$\alpha$ is the simple restricted root of $R_{[L,L],\theta}$, then
$(\omega,\alpha^{\vee})=1$ (see \cite{Ru2}, Theorem 2.2). Let
$\varpi^{\vee}$ be  a generator of $\chi_{*} (Z(L)^{0}/(Z(L)^{0}\cap
H))$; suppose also that there exists
$x':=lim_{t\rightarrow0}\varpi^{\vee}(t)\cdot x_{0}$ in the closure
of $T\cdot x_{0}$ in $Z$.  Then $(\omega,\varpi^{\vee})=1$ because
$(Z(L)^{0}/(Z(L)^{0}\cap H))\cup x'$
is contained in $Z$ (and $Z$ has a unique isotopic component as
$Z(L)^{0}$-representation). To determine the sign of
$\varpi^{\vee}$, we will use the fact that
$(\lambda,\varpi^{\vee})\geq 0$ for every $\lambda\in
\chi(S)\cap\sigma^{\vee}$, where $\sigma$ is the cone associated to
the closure of $T\cdot x_{0}$ in $Z$ and $\sigma^{\vee}$ is its dual
cone. In \cite{Ru2} is proved that $\sigma^{\vee}$ is equal to
$W_{L,\theta}\cdot C(X)^{\vee}$, where $W_{L,\theta}$ is the Weyl
group of the restricted root system of $L$ (see the proof of Lemma
2.8 in \cite{Ru2}).

In the following of this section, let $X$ be a smooth projective
symmetric variety with Picard number one, over which $Aut^{0}(X)$
acts non-transitively. We will prove in \S3 that $X$ has rank 2; for
the moment let us assume it. Suppose also that $Aut^{0}(X)$
stabilizes all the closed $G$-orbits in $X$. We define the canonical
completion of $G/H$ as the simple symmetric variety associated to
the cone $(-C^{+},\emptyset)$. Let $\widetilde{X}$ be the
decoloration of $X$, namely the minimal toroidal variety with a
proper map $\pi:\widetilde{X}\rightarrow X$ which extents the
identity over $G/H$: it is the canonical variety if $X$ is simple
and the blow-up of the canonical variety in the closed orbit
otherwise. Moreover, $\widetilde{X}$ is  the blow-up of $X$ along
the closed orbits (see Theorem 3.3 in \cite{Br2}); in particular,
$\widetilde{X}$ is smooth and the $G$-stable divisor of $X$ is the
image of a $G$-stable divisor of $\widetilde{X}$. Let
$G/\widetilde{P}^{\,\text{-}}$ be a closed orbit of $\widetilde{X}$.

We claim  that, if $Aut^{0}(G/\widetilde{P}^{\,\text{-}})=G/Z(G)$,
then $Lie(Aut^{0}(X))$ is isomorphic to $Lie(G)$. The group
$Aut^{0}(X)$ is isomorphic to $Aut^{0}(\widetilde{X})$. Indeed
$Aut^{0}(X)$ is contained in $Aut^{0}(\widetilde{X})$ because the
closed orbits are stable under the action of $Aut^{0}(X)$. Moreover,
by a result of A. Blanchard, $Aut^{0}(\widetilde{X})$ acts on $X$ in
such a way that the projection is equivariant (see \cite{A2}, $\S
2.4$). Since $\widetilde{X}$ is complete and $Pic\,(\widetilde{X})$
is discrete, the group $Aut^{0}(\widetilde{X})$ is linear algebraic
and its Lie algebra is the space of global vector fields, namely
$H^{0}(X,\mathcal{T}_{\widetilde{X}})$.

We want to prove that $Aut^{0}(\widetilde{X})$ stabilizes all the
$G$-orbits in $\widetilde{X}$. This fact implies that
$Aut^{0}(\widetilde{X})$ is reductive by a result of M. Brion (see
\cite{Br4}, Theorem 4.4.1). First, suppose   $X$ simple; we will
prove in \S\ref{rank2omog} that $H$ is $N_{G}(G^{\theta})$. Thus
$\widetilde{X}$ is a wonderful $G$-variety. Moreover,
$Aut^{0}(\widetilde{X})$ is semisimple, $\widetilde{X}$ is a
wonderful $Aut^{0}(X)$-variety and the set of colors is $D(G/H)$
(see \cite{Br4}, Theorem 2.4.2). In \cite{Br4} it is also determined
the automorphism group of the wonderful completion of a simple
adjoint group $\overline{G}$ (see \cite{Br4}, Example 2.4.5); in
particular, it is proved that $Aut^{0}(\overline{G})$ is
$\overline{G}\times\overline{G}$  if $\overline{G}\neq PSL(2)$.
Coming back to our problem, the wonderful $G$-variety
$\widetilde{X}$ has two stable prime divisors: the exceptional
divisor $E$ and the strict transform $\widetilde{D}$ of the
$G$-stable divisor $D$ of $X$. The divisor $\widetilde{D}$
corresponds to a simple restricted root (see \cite{dCP1}, Lemma
2.2), which is not a dominant weight because the restricted root
system is irreducible of rank 2 (see Theorem~\ref{classif}). Thus
$\widetilde{D}$ is fixed by $Aut^{0}(\widetilde{X})$ (see
\cite{Br4},  Theorem 2.4.1), so $Aut^{0}(\widetilde{X})$ fixes all
the $G$-orbits in $\widetilde{X}$.

We can proceed similarly in the case where  $X$ is non-simple.
Indeed, using \cite{Br1}, Proposition 3.3, \cite{V1}, \S2.4
Proposition 1 and Proposition 2, one can easily prove that
$H^{0}(\widetilde{X},\mathcal{O}(F))=\mathbb{C}s_{F}$ for each
$G$-stable divisor $F$ of $\widetilde{X}$ (here $s_{F}$ is a global
section with divisor $F$). Since $\xi\cdot s_{F}$ is a scalar
multiple of $s_{F}$ for any $\xi\in
H^{0}(\widetilde{X},\mathcal{T}_{\widetilde{X}})=Lie(Aut^{\circ}(\widetilde{X}))$,
$F$ is stabilized by $Aut^{\circ}(\widetilde{X})$ as in the proof of
Theorem 2.4.1 in \cite{Br4} (see also Proposition 4.1.1 in
\cite{BB}). 

Notice that an element of $Z(G)$ acts trivially on $G/H$ (and on
$\widetilde{X}$) if and only if it belongs to $H$, thus
$Aut^{\circ}( \widetilde{X} )$ contains $G/(Z(G)\cap H)$. Assume
that $Aut^{0}(G/\widetilde{P}^{\,\text{-}})$ is $G/Z(G)$, then the
restriction $\psi: Aut^{\circ}(\widetilde{X})\rightarrow
Aut^{\circ}(G/\widetilde{P}^{\,\text{-}})=G/Z(G)$ is an isogeny,
because $Aut^{0}(\widetilde{X})$ stabilizes all the $G$-orbits.
Indeed, $ker\, \psi$ centralizes $G/(Z(G)\cap H)$ (because
$Aut^{0}(\widetilde{X})$ is reductive) and stabilizes $G/H$, so it
is contained in the finite group $Aut_{G}(G/H)\cong N_{G}(H)/H$.
Therefore $Aut^{\circ}( X )$ is $G/(Z(G)\cap H)$.

\rema\label{aut0}  To summarize, given $X$ over which $Aut^{0}(X)$
acts non-transitively, we have to prove:
\begin{itemize}
\item  $Aut^{0}(X)$ stabilizes all the closed $G$-orbits in
$X$;
\item $X$ has rank 2;
\item $Aut(G/\widetilde{P}^{\,\text{-}})=G/Z(G)$ for a closed orbit
$G/\widetilde{P}^{\,\text{-}}$ of $\widetilde{X}$.
\end{itemize}
\emph{If the previous conditions  are verified, then $Aut^{0}(X)$ is
$G/(Z(G)\cap H)$.}

In the following of this section we suppose that the previous
conditions are verified and we try to study the full automorphism
group $Aut(X)$. This group permutes the $Aut^{0}(X)$-orbits; in
particular  it stabilizes the open orbit and the codimension one
orbit. Furthermore, $Aut(X)$ stabilizes $G/H$, thus it is contained
in $Aut(G/H)$. We will prove the converse. Notice that $Aut(X)$ is
contained in $Aut(\widetilde{X})$, because $\widetilde{X}$ is the
blow-up of $X$ along a $Aut(X)$-stable (eventually non-connected)
subvariety.

\rema\label{autXtilta} \emph{The automorphism group of
$\widetilde{X}$ is isomorphic to the automorphism group of $G/H$.}
Indeed, $Aut(\widetilde{X})$ stabilizes $G/H$, so it is contained in
$Aut(G/H)$. Moreover, we can extend every automorphism of $G/H$  to
an automorphism of $\widetilde{X}$, because $\widetilde{X}$ is
normal and  toroidal (these facts imply that  each (non-open)
$G$-orbit $\mathcal{O}$ is contained in the closure of an orbit of
dimension equal to $dim\,\mathcal{O}+1$).   To study $Aut(X)$,  we
have to determine the automorphisms of $\widetilde{X}$ which descend
to $X$.

Now, we study $Aut(\widetilde{X})$ (and $Aut(G/H)$).  For any
$\varphi$ in $Aut(\widetilde{X})$ we define $\tilde{\varphi}\in
Aut_{alg}(Aut^{\circ}(\widetilde{X} ))\equiv Aut_{alg}(G)$ by
\[\tilde{\varphi}(g)\cdot x=\varphi(g\cdot\varphi^{-1}(x)),\ \ \ g\in Aut^{0}(X),\ x\in \widetilde{X}.\]
The kernel of $\varphi\rightarrow\tilde{\varphi}$ is the group of
equivariant automorphism of $\widetilde{X}$ (and $G/H$), thus it is
isomorphic to $N_{G}(H)/H$. We will prove that $N_{G}(H)/H$ is
simple (or trivial); moreover $Z(G/Z(G)\cap H)$ is trivial only if
$N_{G}(H)/H$ is it. Thus $ker(\varphi\rightarrow\tilde{\varphi})$ is
$Z(G/Z(G)\cap H)=Z(G)/(Z(G)\cap H)$.

Let $T''$ be a maximal torus of $H$, then its centralizer
$T':=C_{G}(T'')$  is a maximal torus of $G$ (see \cite{T2}, Lemma
26.2); moreover $T''$ contains a regular one-parameter subgroup
$\lambda$ of $G$. Thus $B':=P(\lambda)$ is a Borel subgroup of $G$
and $B'':=(B'\cap H)^{0}$ is a Borel subgroup of $H$. The group
$Aut_{alg}(G)$ is  generated by $\overline{G}:=G/Z(G)$ and
$E=\{\psi\in Aut_{alg}(G): \psi(T')=T' \textrm{ and }
\psi(B')=B'\}$; the intersection of $E$ with $\overline{G}$ is
$\overline{T}\,':=T'/Z(G)$, so $Aut_{alg}(G)$ is the semidirect
product of $\overline{G}$ and $E':=\{\psi\in E: \psi(t)=t\ \forall
t\in T'\}$. Observe that every $\psi\in E$ induces an automorphism
of the Dynkin diagram (with respect to $T'$ and $B'$), moreover such
automorphism is trivial if and only if $\psi$ belongs to
$\overline{G}$. Furthermore, $Aut(\widetilde{X})$ is generated by
$G/(Z(G)\cap H)$ and by the stabilizer  of $x_{0}$ in
$Aut(\widetilde{X})$. Notice that, given $\varphi\in
Aut_{x_{0}}(\widetilde{X})$, then $\widetilde{\varphi}$ belongs to
$E$, up to composing $\varphi$ with an element of $H$.

\rema \label{aut} Suppose that $Aut_{G}(\widetilde{X})$ is
$Z(G/(Z(G)\cap H))$; in particular it is contained in
$Aut^{0}(\widetilde{X})$. Then $Aut(G/H)$ is generated by
$Aut^{0}(G/H)=G/(Z(G)\cap H)$ and by the subgroup $K=\{\varphi\in
Aut(G/H): \varphi(x_{0})=x_{0} \textrm{ and } \widetilde{\varphi}\in
E\}$. Moreover the map $\varphi\rightarrow\tilde{\varphi}$
restricted to $K$ is injective, because $Z(Aut(G/H))$ does not fix
$x_{0}$. Observe that any automorphism of $G$ stabilizing $H$
induces an automorphism of $G/H$, so $K$ is isomorphic to
$K':=\{\psi\in E:\ \psi(H)=H\}$. Besides $\theta$ belongs to $K$
(see \cite{T2}, Lemma 26.2).

Now, we explain how to prove that every automorphism of
$\widetilde{X}$ descends to an automorphism of $X$. Let $\varphi$ in
$K$, then
$\varphi(g\widetilde{P}')=\widetilde{\varphi}(g)\widetilde{P}'$ for
every $g\widetilde{P}'\in G/\widetilde{P}'$. First, suppose
$\widetilde{X}$ (and $X$) simple; in particular, $\varphi $
stabilizes the closed orbit of $\widetilde{X}$. Let $x_{1}$ be a
point (in the closed orbit of $\widetilde{X}$) fixed by $B'$ and let
$\widetilde{P}'$ be the stabilizer of $x_{1}$ in $G$. Then  $
\varphi $ fixes $x_{1}$, because $\widetilde{\varphi}\in E$. Let
$G/P' $ be the image of $G/\widetilde{P}'\subset \widetilde{X}$ in
$X$; we have to prove that $\widetilde{\varphi}(P')$ is $P'$. Thus,
it is sufficient  that $\varphi$ induces the trivial automorphism of
the Dynkin diagram (with respect to $T'$ and $B'$).

Finally, suppose  that $X$ contains two closed orbits. In this case
the restricted root system has type $A_{2}$, because of Theorem
\ref{classif} and Proposition \ref{A1+A1}. Furthermore, the closed
orbits are $G/P(\omega_{1})$ and $G/P(\omega_{2})$; in particular,
they are not $G$-isomorphic. Let $\mathcal{O}_{1}$ and
$\mathcal{O}_{2}$ be the inverse images in $\widetilde{X}$
respectively of $G/P_{1}$ and $G/P_{2}$; one can easily show that
they are isomorphic. We want to prove that any automorphism in $K$
exchange $\mathcal{O}_{1}$ with $\mathcal{O}_{2}$. Let
$\widetilde{D}$ be the strict trasform of $D$ and let $E_{1}$,
$E_{2}$ be the exceptional divisors of $\widetilde{X}$; we can
suppose that $\mathcal{O}_{i}$ is the intersection of
$\widetilde{D}$ with $E_{i}$. Thus it is sufficient to prove that
$\varphi$ exchange $E_{1}$ with $E_{2}$. For the moment, assume this
fact and     let $x_{i}$ be a point in $\mathcal{O}_{i}$ stabilized
by $B'$, let $P_{i}'$ be the stabilizer of $\pi(x_{i})$ and  let
$P'$ be the stabilizer of $x_{1}$ (and $x_{2}$). We have
$\varphi(x_{1})=x_{2}$, because $\widetilde{\varphi}$ belongs to
$E$.  Now, it is sufficient to prove that $\widetilde{\varphi}$ is
associated to a non-trivial automorphism of the Dynkin diagram.
Indeed, in this case $\widetilde{\varphi}(P'_{1})$ is a parabolic
subgroup of $G$ containing $B'$, distinct by $P'_{1}$ and with the
same dimension of $P'_{1}$; thus it is $P_{2}'$.

\rema \label{autX} Summarizing, to prove that  $Aut(X)$ coincides
with $Aut(G/H)$ (and with $Aut(\widetilde{X})$), it is sufficient
to show that:

\begin{itemize}
\item   $N(H)/H$ is simple or trivial;
\item if $N(H)/H$ is non-trivial, then also $Z(G)/(Z(G)\cap H)$ is
non-trivial;
\item if $X$ is  simple (and non-homogeneous), then $R_{G,\theta}$ has   type
$G_{2}$. Moreover, given $\varphi\in K$ then  $\widetilde{\varphi}$
induces a trivial automorphism of the Dynkin diagram (associated to
$T'$ and $B'$);
\item if $X$ is non-simple and  $\varphi\in K$, then $\varphi$
exchange $E_{1}$ with $E_{2}$. Moreover, $\widetilde{\varphi}$
induces a non-trivial automorphism of the Dynkin diagram (associated
to $T'$ and $B'$).
\end{itemize}
To study $Aut(X)$, we have to determine the group $K\subset
Aut(G/H)$. We will prove that $\theta$ is always contained $Aut
(X)$. If, moreover,  $G/H$ is not isomorphic to $SL_{3}$, then
$K'\cap E'$ is $E'$; thus  the map $\varphi\rightarrow
\widetilde{\varphi}$ is surjective (if $G/H\neq SL_{3}$). In
following subsections we will study $Aut(X)$ by a case-to-case
analysis.

\subsection{Homogeneous varieties}\label{rank2omog}

We begin studying the symmetric varieties of rank 2 which are
homogeneous.

\begin{prop}\label{AIII}
The smooth completion of $SL_{4}/N_{SL_{4}}(S(L_{2}\times L_{2}))$
with Picard number one is isomorphic to $\mathbb{G}_{2}(6)$.
\end{prop}


{\em Proof.} The symmetric varieties
$SL_{4}/N_{SL_{4}}(S(L_{2}\times L_{2}))$ has type $AIII$. The group
$SL_{4}$ acts on the six-dimensional space
$\bigwedge^{2}(\mathbb{C}^{4})$, thus it acts on
$\mathbb{G}_{2}(\bigwedge^{2}\mathbb{C}^{4})$. The stabilizer of the
space generated by $e_{1}\wedge e_{2}$ and $e_{2}\wedge e_{3}$ is
$N_{SL_{4}}(S(L_{2}\times L_{2}))$, thus
$SL_{4}/N_{SL_{4}}(S(L_{2}\times L_{2}))$ is contained in
$\mathbb{G}_{2}(\bigwedge^{2}\mathbb{C}^{4})$. Moreover
$SL_{4}/N_{SL_{4}}(S(L_{2}\times L_{2}))$ has the same dimension of
$\mathbb{G}_{2}(\bigwedge^{2}\mathbb{C}^{4})$, so the Grasmannian is
the unique smooth completion of $SL_{4}/N_{SL_{4}}(S(L_{2}\times
L_{2}))$ with Picard number one. $\square$

\begin{prop}\label{EIV}
The smooth completion of $E_{6}/N_{E_{6}}(F_{4})$ with Picard number
one is isomorphic to $\mathbb{P}(\mathbb{J}_{3}(\mathbb{O}))$.
\end{prop}

{\em Proof.}  The symmetric varieties $E_{6}/N_{E_{6}}(F_{4})$ has
type $EIV$ and $\mathbb{J}_{3}(\mathbb{O})$ is the 27-dimensional
irreducible representation of $E_{6}$ corresponding to the first
fundamental weight. The subgroup $F_{4}$ of $E_{6}$ is isomorphic to
the group of automorphism of $\mathbb{J}_{3}(\mathbb{O})$; in
particular $F_{4}$ fixes the identity matrix. Thus
$\mathbb{P}(\mathbb{J}_{3}(\mathbb{O})) $ contains the
26-dimensional variety $G/H$. $\square$

\begin{prop}\label{CIIii-agg}
The smooth completion of $Sp_{8}/N_{Sp_{8}}(Sp_{4}\times Sp_{4})$
with Picard number one is isomorphic to
$E_{6}/P_{1}\equiv\mathbb{P}^{2}(\mathbb{O})$.
\end{prop}

{\em Proof.}  The symmetric variety $Sp_{8}/N_{Sp_{8}}(Sp_{4}\times
Sp_{4})$ has type  $CII$. Let $(G',\sigma)$ be an involution of type
$EI$, where $G'$ is the simply connected simple  group of type
$E_{6}$ and $G'^{\sigma}$ has type $C_{4}$. Choose a maximally
$\sigma$-split torus and  a Borel subgroup of $G'$ as in
\S˜\ref{1def}; then  $\sigma$ acts as $-id$ over $R_{G'}$ and the
parabolic subgroup $P:=P(\varpi_{1})$ of $G'$ is opposed to
$\sigma(P)$. Furthermore, $P\cap \sigma(P)$ is a Levi subgroup of
$P$ containing $P^{\sigma}$; the derived subgroup of $P\cap
\sigma(P)$ has type $D_{5}$ while $(P^{\sigma})^{0}$  has type
$B_{2}\times B_{2}$. Hence $G'^{\sigma}/P^{\sigma}$ is isomorphic to
$Sp_{8}/N_{Sp_{8}}(Sp_{4}\times Sp_{4})$, up to quotient by a finite
group. Notice that $G'/P$ is a smooth completion of
$G'^{\sigma}/P^{\sigma}$ with Picard number one. The variety
$G'^{\sigma}/P^{\sigma}$ cannot be isomorphic to
$Sp_{8}/(Sp_{4}\times Sp_{4})$ because the unique smooth completion
of $Sp_{8}/(Sp_{4}\times Sp_{4})$ with Picard number one is
isomorphic to $\mathbb{G}_{4}(8)$ (see Proposition \ref{CII} in
\S~\ref{cap BCl+Cl} and Theorem \ref{classif}). $\square$

Now, we study the symmetric varieties whose restricted root system
is reducible. If the restricted root system of $G/H$ has type
$A_{1}$, then $G/N_{G}(G^{\theta})$ is isomorphic to
$SO_{n}/N_{SO_{n}}(SO_{n-1})$ for an appropriate $n\geq 3$. If
$n=3$,  the type of $(G,\theta)$ is $AI$; if $n=4$  then
$G/N_{G}(G^{\theta})$ is isomorphic to $PSL_{2}$; if $n=6$,  the
type of $(G,\theta)$ is $AII$; if $n=2k$ with $k\geq4$,  the type of
$(G,\theta)$ is $DII$;  if $n=2k+1$ with $k\geq2$, the type of
$(G,\theta)$ is $BII$.

\begin{prop}\label{A1+A1}
If $R_{G,\theta}$ has type $A_{1}\times A_{1}$,   then
$G/N_{G}(G^{\theta})$ is isomorphic to
$SO_{n}/N_{SO_{n}}(SO_{n-1})\times SO_{m}/N_{SO_{m}}(SO_{m-1})$ and
the smooth projective embedding of $G/H$ with Picard number one is
isomorphic to $\mathbb{I}\mathbb{G}_{1}(n+m)$.
\end{prop}

First we want to describe $H$: write
$(G,\theta)=(G_{1},\theta)\times (G_{2},\theta)$ and let $g_{i}\in
N_{G_{i}}(G_{i}^{\theta})$ be a representant of the non-trivial
element of $N_{G_{i}}(G_{i}^{\theta})/G_{i}^{\theta}$ for each $i$.
Then $H$ is generated by $G^{\theta}$ and $(g_{1},g_{2})$.

{\em Proof.} The group $SO_{n}\times SO_{m}$ is contained in
$SO_{n+m}$. Let $\{e_{1},...,e_{n+m}\}$ be an orthonormal basis of
$\mathbb{C}^{n+m}$ such that $SO_{n}\subset
GL(span_{\mathbb{C}}\{e_{1},...,e_{n}\})$ and $SO_{m}\subset
GL(span_{\mathbb{C}}\{e_{n+1},...,e_{n+m}\})$. The connected
component   of the stabilizer of $e_{1}+ie_{n+1}$ in $SO_{n}\times
SO_{m}$ is $SO_{n-1}\times SO_{m-1}$ (see also \cite{dCP1}, Lemma
1.7). Furthermore, $SO(\mathbb{C}^{n+m})\cdot
[e_{1}+ie_{n+1}]\subset \mathbb{P}( \mathbb{C}^{n+m})$ is
$\mathbb{I}\mathbb{G}_{1}(n+m)$ because $e_{1}+ie_{n+1}$ is an
anisotropic vector. (Notice that if $R_{G,\theta}$ has type
$A_{1}\times A_{1}$, then there is a unique subgroup
$G^{\theta}\subset H\subset N_{G}(G^{\theta})$ such that $G/H$ has a
smooth completion with Picard number one). $\square$

\subsection{Restricted root system of type $G_{2}$}\label{cap-g2}

Now, we study the symmetric varieties whose restricted root system
has type $G_{2}$. First of all, we determine their connected
automorphism group; in particular, we prove that such varieties are
non-homogeneous.

\begin{lem}\label{lemG} Suppose that $R_{G,\theta}$ has type
$G_{2}$. Then the connected automorphism group $Aut^{0}(X)$ is
isomorphic to $G$.
\end{lem}

{\em Proof of Lemma~\ref{lemG}.} First, we prove that
$Aut^{\circ}(X)$ does not act transitively over $X$. This simple
variety is  associated to the colored cone
$(\sigma(\alpha_{2}^{\vee}, -\omega_{2}^{\vee}),$ $
\{D_{\alpha_{2}^{\vee}}\})$  and its closed orbit is isomorphic to
$G/P(-\omega_{1})$. Let $\omega$ be the highest weight of $Z$ (with
respect to $L$). The simple restricted root of $R_{[L,L],\theta}$ is
$\alpha_{2}$, so $(\omega,\alpha_{2}^{\vee})=1$ (see
\S\ref{general}). Moreover, $Z(L)^{0}/(Z(L)^{0}\cap H)$ is the
one-parameter subgroup of $T/T\cap H$ corresponding to
$-\omega_{1}^{\vee}$; indeed $(-\omega_{1} , \omega_{1}^{\vee})<0 $
and $-\omega_{1} $ belongs to the dual cone $\sigma^{\vee}$
associated to the closure of $T\cdot x_{0}$ in $Z$. Since
$Z(L)^{0}/(Z(L)^{0}\cap H)\equiv Z(L)^{0}\cdot x_{0}$ is contained
in $Z$, we have $1=(\omega,-\omega_{1}^{\vee})=
(a\omega_{1}+\omega_{2},-2\alpha_{1}-\alpha_{2})=-2a-1$, so $\omega$
is $-\omega_{1}+\omega_{2}$. Thus $Aut^{0}(X)$ is isomorphic to $G$
by the Remark \ref{aut0}.  Observe that the connected automorphism
group of the closed orbit  of $\widetilde{X}$ is $G$, while the
connected automorphism group of the closed orbit of $X$ is either
$SO_{7}$ or $SO_{7}\times SO_{7}$. $\square$

\begin{prop}\label{G}
Suppose that  $(G,\theta)$ has type $G$ (and that $G/H$ is
isomorphic to $G_{2}/(SL_{2}\times SL_{2})$). Then $Aut(X)$ has
three orbits in $X$, is connected and is isomorphic to $G_{2}$.
Furthermore, $X$ is intersection of hyperplane sections of
$\mathbb{G}_{3}(7)$. More precisely, $X$ is the intersection in
$\mathbb{P}(\bigwedge^{3}V(\omega_{1}))$ of $\mathbb{G}_{3}(7)$ with
$\mathbb{P}(V(\omega_{1}))$. Moreover, $X$ is the subvariety of
$\mathbb{G}_{3}(Im\, \mathbb{O})$ 
parametrizing the subspaces $W$ such that $\mathbb{C}1\oplus W$ is a
subalgebra of $\mathbb{O}$ 
isomorphic to $\mathbb{H}$. 

\end{prop}

{\em Proof of Proposition~\ref{G}.} The center of $G_{2}$ is
trivial, thus $Aut(X)$ is contained in $Aut_{alg}(G_{2})$ (see the
observations after the Remark \ref{autXtilta}); such group is
connected, so also $Aut(X)$ is connected.

Now, we prove that $X$ is "contained'' in $\mathbb{G}_{3}(7)$. There
is an involution of $SO_{7}$ that extends $\theta$; again, we denote
it by $\theta$. We have $SO_{7}^{\theta}=S(O_{3}\times SO_{4})$,
thus  $G_{2}/(SL_{2}\times SL_{2})$ is a closed subvariety of
$SO_{7}/N_{SO_{7}}(S(O_{3}\times SO_{4}))$. There is a unique smooth
completion of $SO_{7}/N_{SO_{7}}(S(O_{3}\times SO_{4}))$ with Picard
number one and is isomorphic to $\mathbb{G}_{3}(7)$ (see Theorem
\ref{classif} and Proposition \ref{BI+DI}). We have to introduce
some notations. Let $V$ be a 7-dimensional vector space and let
$\{e_{-3},e_{-2},e_{-1},e_{0},e_{1},e_{2},e_{3}\}$ be a basis of
$V$. Let $q$ be the symmetric bilinear form associated to the
quadratic form
$(e^{*}_{0})^{2}+e^{*}_{1}e^{*}_{-1}+e^{*}_{2}e^{*}_{-2}+e^{*}_{3}e^{*}_{-3}$
and let $\varpi$ be the trilinear form $e^{*}_{0}\wedge
e^{*}_{1}\wedge e^{*}_{-1}+e^{*}_{0}\wedge e^{*}_{2}\wedge
e^{*}_{-2}+e^{*}_{0}\wedge e^{*}_{3}\wedge
e^{*}_{-3}+2\,e^{*}_{1}\wedge e^{*}_{2}\wedge
e^{*}_{3}+2\,e^{*}_{-1}\wedge e^{*}_{-2}\wedge e^{*}_{-3}\in
\bigwedge^{3}V^{*}$. The subgroup $G$ of $SL(V)$ composed by the
linear transformations which preserve $q$ and $\varpi$ is the simple
group of type $G_{2}$; moreover, we can realize $SO_{7}$ as
$SO(V,q)$. The vector space $V$ is  the standard representation
$V(\omega_{1})$ of $G$ and we can suppose that
$\{e_{-3},e_{-2},e_{1},e_{0},e_{1},e_{2},e_{3}\}$ is a basis of
weight  vector for an appropriate maximal torus $T$ of $G$.
Moreover, we can choose a Borel subgroup $B$ of $G$ so that the
weight  of $e_{i}$ is a positive root (respectively is 0) if and
only if $i>0$ (respectively $i=0$). The Grassmannian
$\mathbb{G}_{3}(7)$ is contained in $\mathbb{P}(\bigwedge^{3}V)$ and
$\bigwedge^{3}V$ is isomorphic to $V\oplus
V(2\omega_{1})\oplus\mathbb{C}$ as $G_{2}$-representation. The
subrepresentation of $\bigwedge^{3}V$ isomorphic to $V$ has the
following basis of $T$-weights: $\{ 
e_{-2}\wedge e_{-3}\wedge e_{0}-e_{1}\wedge e_{2}\wedge
e_{-2}-e_{1}\wedge e_{3}\wedge e_{-3},
e_{2}\wedge e_{-2}\wedge e_{-3}-e_{1}\wedge e_{2}\wedge
e_{0}+e_{1}\wedge e_{-1}\wedge e_{-3},
e_{3}\wedge e_{-2}\wedge e_{-3}-e_{1}\wedge e_{3}\wedge
e_{0}-e_{1}\wedge e_{-1}\wedge e_{-2},
e_{1}\wedge e_{2}\wedge e_{3}-e_{-1}\wedge e_{-2}\wedge e_{-3},
e_{2}\wedge e_{3}\wedge e_{-3}-e_{1}\wedge e_{2}\wedge
e_{-1}+e_{0}\wedge e_{-1}\wedge e_{-3},
-e_{2}\wedge e_{3}\wedge e_{-2}-e_{1}\wedge e_{3}\wedge
e_{-1}-e_{0}\wedge e_{-1}\wedge e_{-2},
-e_{1}\wedge e_{-1}\wedge e_{-2}-e_{2}\wedge e_{3}\wedge
e_{0}-e_{3}\wedge e_{-1}\wedge e_{-3} \}$. Let $X''$ be the closure
of $G_{2}/(SL_{2}\times SL_{2})$ in $\mathbb{G}_{3}(7)$ and let $X'$
be the intersection of $\mathbb{G}_{3}(7)$ with
$\mathbb{P}(V(2\omega_{1})\oplus\mathbb{C})$; observe that $X''$ is
contained in $X'$ because $V=V(\omega_{1})$ is not spherical (see
\cite{T2}, Proposition 26.4). We claim that $X$ is the normalization
of $X''$. Indeed $\mathbb{P}(V(2\omega_{1})\oplus\mathbb{C})$ has
two closed $G_{2}$-orbits: one is isomorphic to
$G_{2}/P(\omega_{1})$ and the other one is the point
$\mathbb{P}(\mathbb{C})$. Therefore $X'$ contains one closed
$G_{2}$-orbit, otherwise $\mathbb{G}_{3}(7)$ would contain a
$G_{2}$-fixed point, in particular $V$ would  be reducible as
$G_{2}$-representation. Thus, the normalization of $X''$ is a simple
variety with closed orbit $G/P(\omega_{1})$, so it is $X$.

We want to prove that $X''$ coincides with $X'$ and is smooth (so it
coincides also with $X$). Notice that $X'$ is connected, because it
is $G_{2}$-stable and contains a unique closed $G_{2}$-orbit. Hence,
it is sufficient to prove that $X'$ has dimension 8 and is smooth in
a neighborhood of a point belonging to $G_{2}/P(\omega_{1})$. One
can verify that $e_{1}\wedge e_{-2}\wedge e_{-3}$ is a highest
weight vector of $V(2\omega_{2})\subset \bigwedge^{3}V$. Let $A$ be
the affine open subset of $\mathbb{G}_{3}(7)$ composed by the
subspaces with basis $\{e_{j}+\sum a_{i,j}e_{i}\}_{j=1,-2,-3;
i=2,3,0,-1}$, where $(a_{i,j})$ varies  in $\mathbb{C}^{12}$. The
subvariety $A\cap\mathbb{P}(V(2\omega_{1}\oplus\mathbb{C}))=A\cap
X'$ of $A$ has equations:
\begin{itemize}
\item[ ] $a_{1,0}=-a_{3,2}+a_{2,3},$
\item[ ] $a_{2,-1}=-a_{1,2}+T_{(2,3),(2,0)},$
\item[ ] $a_{3,-1}=-a_{1,3}+T_{(2,3),(3,0)},$
\item[ ] $a_{1,-1}=T_{(2,3),(2,3)},$
\item[ ] $T_{(1,2),(2,3)}-T_{(2,3),(2,-1)}+T_{(1,2),(0,-1)},$
\item[ ] $T_{(1,3),(2,3)}-T_{(2,3),(3,-1)}+T_{(1,3),(0,-1)},$
\item[ ] $a_{3,-1}+T_{(1,2,3),(2,3,0)}+T_{(1,2),(3,-1)}$
\end{itemize} (where $T_{(h,k),(n,m)}$ is the minor of
$(a_{i,j})$  extracted by the $h$-th and $k$-th row and by $n$-th
and $m$-th column). The closed subset $A'$ of $A$  defined by the
the first four equation is the graph of a polynomial map, thus it is
smooth of dimension 8. Hence the last three equation are identically
verified on $A'$, because $X'' $ has dimension 8 (and $A\cap
X''\subset A\cap X'\subset A'$).  Therefore $A\cap X'$ coincides
with $A'$ and is smooth.

Now, we prove the last statement of the proposition. Identify $V$
with $Im (\mathbb{O})$ 
and define the associator
$[\cdot,\cdot,\cdot]:\bigwedge^{3}(\mathbb{O})
\rightarrow (\mathbb{O})$ 
as the linear map such that $[a,b,c]=(ab)c-a(bc)$. This map is
$G_{2}$-equivariant and has kernel $V(2\omega_{1})\oplus\mathbb{C}$.
Thus $X$ parametrizes the 3-dimensional subspaces $W$ of $Im
(\mathbb{O})$ 
over which $[\cdot,\cdot,\cdot]$ is zero. Furthermore,
$[1,\cdot,\cdot]$ is identically zero. Let $W$ be a subspace
associated to a point in $X$ and let $\overline{W}$ be the
subalgebra of $\mathbb{O}$ 
generated by $W$ and $1$. It can be either the entire algebra
$\mathbb{O}$ 
or a subalgebra of dimension four. But, if it is the whole algebra,
then $\mathbb{O}$ 
is generated by four elements which associates between them; a
contradiction (recall that $\mathbb{O}$
and $\overline{W}$ are composition algebras).   Thus, $\overline{W}$
is a composition algebra of dimension four, so it is isomorphic to
$\mathbb{H}$. 
$\square$

Let $(V_{1},q_{1})$ and  $(V_{2},q_{2})$  be  copies respectively of
$(V,q)$ and $(V,-q)$, where $V$ is as in the previous proposition;
we can suppose $G_{2}\times G_{2}\subset SO(V_{1})\times SO(V_{2})$.
Moreover, let $W_{i}$ be a maximal anisotropic subspace of $V_{i}$
for both the $i$. Let $W$ be a maximal anisotropic subspace of
$V_{1}\oplus V_{2}$ which contains $W_{1}\oplus W_{2}$.

\begin{prop}\label{G2}
Suppose that $G/H$ is is isomorphic to the simple group of type
$G_{2}$. Then  $Aut(X)^{0}$ is $G_{2}\times G_{2}$, while $Aut(X)$
is generated by $Aut(X)^{0}$ together with the flip
$(g,h)\rightarrow (h,g)$. In particular, $Aut(X)^{0}$ has index two
in $Aut(X)$. Furthermore, $X$ is intersection of  hyperplane
sections of $\mathbb{I}\mathbb{G}_{7}(14)$: more precisely, $X$ is
the intersection in $\mathbb{P}(\bigwedge^{even}(V_{1}\oplus
V_{2}))) \cong\mathbb{P}((V_{1}\otimes V_{2})\oplus V_{1}\oplus
V_{2}\oplus \mathbb{C})$ of $\mathbb{I}\mathbb{G}_{7}(14)$ with
$\mathbb{P}((V_{1}\otimes V_{2})\oplus \mathbb{C})$
\end{prop}

{\em Proof  of Proposition~\ref{G2}}. The group $G$ has type
$G_{2}\times G_{2}$,  $G^{\theta}$ is the diagonal and $X$ has
dimension 14. The automorphism group of $X$ is determined in
\cite{Br4}, Example 2.4.5. Alternatively, one can study it in a very
similar way to the Proposition \ref{G}.

Clearly the involution of $G\times G$ can be  extended to an
involution of $SO(V_{1})\times SO(V_{2})$ which we denote again by
$\theta$; in particular we have $G\cong G\times G/(G\times
G)^{\theta}\subset SO(V_{1})\times SO(V_{2})/(SO(V_{1})\times
SO(V_{2}))^{\theta}\cong SO(V,q )$.

Let $X''$ be the closure of  $G$ in the unique smooth completion of
$SO(V,q)$ with Picard number one.
This last variety is isomorphic to the spinorial variety
$\mathbb{S}_{7}$ and  the application $\Phi:SO(V,q)\hookrightarrow
\mathbb{S}_{7}$ sends an element $g\in SO(V,q)$ to the graph
$graph(g):=\{(v,gv): v\in V\}\subset V\oplus V\equiv V_{1}\oplus
V_{2}$ of $g$ (see also Proposition \ref{Bl+Dl}). We claim that $X$
is the normalization of $X''$ (one could show that $X''$ is normal
by \cite{T1}, Proposition 9).

Let $\varphi:\mathbb{S}_{7}\rightarrow\mathbb{P}(\bigwedge^{even}W)$
be the $Spin_{14}$-equivariant embedding of $\mathbb{S}_{7}$ in the
projectivitation of  a  half-spin representation of $Spin_{14}$.
Write $V_{i}=W_{i}\oplus\widetilde{W}_{i}\oplus\mathbb{C}_{i}$,
$V=V_{1}\oplus V_{2}=W\oplus\widetilde{ W}$ where $W_{i}$,
$\widetilde{W}_{i}$, $\widetilde{W}$ are maximal anisotropic
subspaces such that
$\widetilde{W}_{1}\oplus \widetilde{W}_{2}\subset \widetilde{W}$.
The representation $\bigwedge^{even}W$ is isomorphic to
$\bigwedge^{\bullet} W_{1}\otimes \bigwedge^{\bullet} W_{2}\cong
\bigwedge^{\bullet} W_{1}\otimes (\bigwedge^{\bullet} W_{2})^{*}$ as
$(Spin_{7}\times Spin_{7})$-representation (see the highest weights
and the dimensions). Moreover $\bigwedge^{\bullet}W_{i}$ is
isomorphic to $V_{i}\oplus \mathbb{C}_{i}$ as
$G_{2}$-representation, so $\bigwedge^{\bullet} W_{1}\otimes
\bigwedge^{\bullet} W_{2}$ is isomorphic to $(V_{1}\otimes
V_{2})\oplus V_{1}\oplus V_{2}\oplus \mathbb{C}$ as
$G$-representation. Let $\mathbb{P}$ be the projective subspace of
$\mathbb{P}(\bigwedge^{even}W)$ isomorphic to
$\mathbb{P}((V_{1}\otimes V_{2})\oplus \mathbb{C})$. Observe that
$X''$ is contained in $X':=\mathbb{S}_{7}\cap\, \mathbb{P}\subset
\mathbb{P}(\bigwedge^{even}W)$ because $V_{1}\oplus V_{2}$ does not
contain a line fixed by $G^{\theta}$ (see \cite{T2}, Proposition
26.4).

Observe that $X'$ contains one closed $G$-orbit. Indeed $\mathbb{P}$
contains two closed $G$-orbits: one isomorphic to $G/P(\omega_{1})$
and the other one isomorphic to the $G$-stable point $\mathbb{P}(
\mathbb{C})$. On the other hand, there is not a $G$-stable maximal
isotropic subspace of $V$, so $\mathbb{P}( \mathbb{C})$ is not
contained in $\mathbb{S}_{7}$. Thus $X'$ is connected. Moreover, the
normalization of $X''$ is the simple symmetric variety with closed
orbit $G/P(\omega_{1})$, so it is $X$. We want to prove that $X''$
is smooth and coincides with $X'$; it is   sufficient to prove that
$X'$ is smooth of dimension 12; in this case $X'$ is irreducible, so
it coincides with $X''$ (and $X$). Moreover, it is sufficient to
study $X'$ in a open neighborhood of an arbitrarily fixed point of
$G/P(\omega_{1})$, for example $x=[e_{1}\wedge e_{-2}\wedge
e_{-3}\wedge e_{0}+f_{0}\wedge f_{2}\wedge f_{3}\wedge f_{-1}]$.

Let $\{e_{-3},e_{-2},e_{-1},e_{0},e_{1},e_{2},e_{3}\}$ be a basis of
$V_{1}$ as before and let $\{f_{-3},f_{-2},f_{1},$ $f_{0},$
$f_{1},f_{2},f_{3}\}$ be the corresponding basis of $V_{2}$. Let $u$
be $e_{0}+f_{0}$; we can suppose that $W$ is generated by
$e_{1},e_{2},e_{3},f_{1},f_{2},f_{3}$ and $u$. The trivial
subrepresentation $\mathbb{C}_{1}$ of $ \bigwedge^{\bullet}W_{1}$ is
spanned by $2\sqrt{-2}1_{W_{i}}+e_{1}\wedge e_{2}\wedge e_{3}$;
moreover, $ W_{1}\oplus \bigwedge^{2}W_{1}$ is contained in the
$G_{2}$-stable   subspace of $ \bigwedge^{\bullet}W_{1}$ isomorphic
to $V_{1}$. An open neighborhood of $x$ in $\mathbb{S}_{7}$ is given
by $U^{\text{-}}\cdot x$, where $U^{\text{-}}$ is the unipotent
radical of standard parabolic subgroup opposite to
$Stab_{Spin_{14}}(x)$ (notice that, as algebraic variety,
$U^{\text{-}}\cdot x$ is isomorphic to $Lie(U^{\text{-}})$ by the
exponential map). The coordinates of $exp(p)\cdot x$ are the
pfaffians of the diagonal minors of $p$. Let $x_{i,j}$ be the
coordinates of the space $M_{14}$ of matrices of order 14 with
respect to the basis $\{e_{1},e_{2},e_{3},u,f_{1},f_{2},f_{3},$
$e_{-1},e_{-2},e_{-3},\frac{1}{2}(e_{0}-f_{0}),f_{-1},f_{-2},f_{-3}\}$
(notice that $ Lie (Spin_{14})\subset M_{14}$). We claim that a open
neighborhood of $x$ in $X'\,\cap\, (U^{\text{-}}\cdot x)$ is the
graph of a polynomial map. Given an skew-symmetric matrix of order
$2n$, let $[i_{1},i_{2},...,i_{2k}]$ be the Pfaffian of the
principal minor extracted from the rows and the columns of indices
$i_{1}<i_{2}<...<i_{2k}$.

The three equations corresponding  to the vectors
$2\sqrt{-2}e_{i}\wedge e_{j}+e_{i}\wedge e_{j}\wedge u\wedge
f_{1}\wedge f_{2}\wedge f_{3}$ ($\in V_{1}$) are, respectively,
$0=x_{i,j}-[i,j,4,5,6,7]=x_{i,j}-(x_{4,7}p_{i,j}+q_{i,j})$, where
the $p_{i,j}$, $q_{i,j}$ are homogeneous polynomials in the
$x_{h,k}$ such that : 1) $(h,k)\neq(4,7)$; 2) either $h>3$ or $k>3$.
Finally, we consider the equation
$0=x_{4,7}-[1,2,3,4]=x_{4,7}-(x_{1,2}x_{3,4}-x_{1,3}x_{2,4}+x_{2,3}x_{1,4})$
associated to $f_{4}\wedge u+f_{4}\wedge e_{1}\wedge e_{2}\wedge
e_{3}$ ($\in V_{2}$). Substituting the first three equation in the
last one, we obtain
$x_{4,7}(1+x_{1,2}f_{1,2}+x_{1,3}f_{1,3}+x_{2,3}f_{2,3})=g$ where
the $f_{i,j}$ and $g$ are polynomial in the $x_{4,j}$ with $j<4$.

Therefore, there is  an open neighborhood $A$ of $x$ in
$U^{\,\text{-}}\cdot x$ where the previous four equation became
$x_{1,2}=h_{1,2}$, $x_{1,3}=h_{1,3}$, $x_{2,3}=h_{2,3}$ and
$x_{4,7}=h_{4,7}$ (here the $h_{i,j}$ are polynomials in the
coordinates different from $x_{1,2}$, $x_{1,3}$, $x_{2,3}$ and
$x_{4,7}$). Observe that the previous equation are independent;  on
the other hand $X'$ has dimension at least 12, because it contains
$X''$. Therefore, the subvariety $A'$ of $A$ obtained imposing the
previous four equation    is smooth with dimension 12, thus it is
equal to $A\cap X''$. 
Hence $X'=X''(=X)$. $\square$

\subsection{Restricted root system of type $A_{2}$ (with
$H=G^{\theta}$)}\label{a2}

Now, we consider the symmetric varieties such that: 1) the
restricted root system has type $A_{2}$; 2) $H=G^{\theta}$. We prove
that they are hyperplane sections of the varieties of the
subexceptional serie of subadjoint varieties (these last varieties
constitute the third line of the geometric version of Freudenthal's
magic square and  are Legendrian varieties).  For the completion of
$SL_{2}$ this result is due to Buczy\'{n}ski (see \cite{Bu}). These
varieties are contained nested in each other. First, we prove that
such varieties are hyperplane sections of  Legendrian varieties.
Then, we study their connected automorphism group; in particular we
show that these varieties are not homogeneous. Finally, we study
their automorphism group.

Recall that $\mathcal{J}_{3}(\mathbb{A})$ is the space of Hermitian
matrices of order three over the complex composition algebra
$\mathbb{A}$. Moreover  $SL_{3}(\mathbb{A})$ is the subgroup of
$GL_{\mathbb{C}}(\mathcal{J}_{3}(\mathbb{A}))$  of complex linear
transformations preserving the determinant, while
$SO_{3}(\mathbb{A})$ is the subgroup of complex linear
transformations preserving the Jordan multiplication. The space
$\mathcal{Z}_{2}(\mathbb{A}):=\mathbb{C}\oplus
\mathcal{J}_{3}(\mathbb{A})\oplus\mathcal{J}_{3}(\mathbb{A})^{*}\oplus\mathbb{C}^{*}$
is an irreducible $Sp_{6}(\mathbb{A})$
-representation. The
closed $Sp_{6}(\mathbb{A})$-orbit in
$\mathbb{P}(\mathcal{Z}_{2}(\mathbb{A}))$ is
$\mathbb{L}\mathbb{G}(\mathbb{A}^{3},\mathbb{A}^{6})$ and  is the
image of the $Sp_{6}(\mathbb{A})$-equivariant rational map:

\[\begin{matrix}
\phi:\mathbb{P}(\mathbb{C}\oplus\mathcal{J}_{3}(\mathbb{A}))&\dashrightarrow&
 \mathbb{P}(\mathcal{Z}_{2}(\mathbb{A}))\\
 (x,P)&\rightarrow &(x^{3},x^{2}P,xcom\,P,det\,P).\end{matrix}\]

The quotient $SL_{3}(\mathbb{A})/SO_{3}(\mathbb{A})$ is a symmetric
variety isomorphic to the image of $SL_{3}(\mathbb{A})\cdot[1,I]$ by
$\phi$; in particular it is contained in the hyperplane section
$X':=\{[x_{1},x_{2},x_{3},x_{4}]\in
\mathbb{L}\mathbb{G}(\mathbb{A}^{3},\mathbb{A}^{6})\subset\mathbb{P}(\mathbb{C}\oplus
\mathcal{J}_{3}(\mathbb{A})\oplus\mathcal{J}_{3}(\mathbb{A})^{*}\oplus\mathbb{C}^{*}):
\ x_{1}=x_{4}\}$. Such section is irreducible because it is the
image by $\phi$ of $\{[x,P]:x^{3}-det(P)\}$; moreover, it has the
same dimension of $SL_{3}(\mathbb{A})/SO_{3}(\mathbb{A})$, thus it
is the closure of $SL_{3}(\mathbb{A})/SO_{3}(\mathbb{A})$ in
$\mathbb{P}(\mathcal{Z}_{2}(\mathbb{A}))$.

\begin{lem}\label{legendr}
The variety $X'$ is isomorphic to the unique smooth completion $X$
of $SL_{3}(\mathbb{A})/SO_{3}(\mathbb{A})$ with Picard number one.
\end{lem}

{\em Proof of Lemma \ref{legendr}}. The variety $X'$ and its
normalization $X''$ have the same number of orbits by \cite{T1},
Proposition 1; moreover, each orbit of  $X''$ have the same
dimension of its image in $X'$. Observe that no orbits of $X''$ have
dimension 0; thus $X'$ has two closed orbits, namely $G/P(\omega_{1}
)$ and $G/P( \omega_{2})$.
One can easily show that $X''$ has a unique orbit $\mathcal{O}$ of
codimension one (otherwise,  by the theory of spherical embeddings,
$X'$ would contain a closed orbit isomorphic to
$G/P(\omega_{1}+\omega_{2})$). Thus $X''=G/H\,\cup \, \mathcal{O}\,
\cup\, G/P(\omega_{1})\,\cup\, G/P(\omega_{2})$.
%
By Proposition 8.2 in~\cite{LM}, we know the possible dimensions for
the singular locus of $X'$. Because it is
$SL_{3}(\mathbb{A})$-stable, one can easily see that $X'$ is smooth
(and coincides with $X''$). Studying the colored fan of $X''$, one
can show easily that $X''$ has Picard number one, so  it is
isomorphic to $X$.
%
%
%
%
%
$\square$

Observe that we have proved the following proposition.

\begin{prop}
The smooth completion of $SL_{3}(\mathbb{A})/SO_{3}(\mathbb{A})$
with Picard number one is contained in the smooth completion of
$SL_{3}(\mathbb{A'})/SO_{3}(\mathbb{A'})$ with Picard number one and
$\mathbb{A'}\subset \mathbb{A}$.  Moreover, we have a commutative
diagram:

\[ \xymatrix{ \overline{SL_{3}/SO_{3}}\ \ar@{^{(}->}[r]\ar@{^{(}->}[d]& \overline{SL_{3}}
\ \ar@{^{(}->}[r]\ar@{^{(}->}[d]& \overline{SL_{6}/Sp_{6}}\
\ar@{^{(}->}[r]\ar@{^{(}->}[d]&
\overline{E_{6}/F_{4}}\ar@{^{(}->}[d]\\\mathbb{LG}_{3}(6) \
\ar@{^{(}->}[r] & \mathbb{G}_{3}(6)\ \ar@{^{(}->}[r]&
\mathbb{S}_{12}\ \ar@{^{(}->}[r]& \mathbb{G}_{3}(\mathbb{O}^{6})\\
}\]
\end{prop}

\begin{lem}\label{lemA2} Suppose that $R_{G,\theta}$ has type
$A_{2}$ and that $H=G^{\theta}$. Then the connected automorphism
group $Aut^{0}(X)$ is isomorphic to $G$, up to isogeny.
\end{lem}

{\em Proof of Lemma~\ref{lemA2}.}  The character group of $\chi(S)$
is  the lattice generated by spherical weights. Let $X'$ be  the
simple variety corresponding to the colored cone
$(\sigma(\alpha_{1}^{\vee}, -\omega_{1}^{\vee}-\omega_{2}^{\vee}),$
$ \{D_{\alpha_{1}^{\vee}}\})$; its  closed orbit  is isomorphic to
$G/P(\omega_{2})$. Let $Z$ be as $\S \ref{general}$ and let
$\omega=a_{1}\omega_{1}+a_{2}\omega_{2}$ be the highest weight of
$Z$ (with respect to $L$). The simple restricted root of
$R_{[L,L],\theta}$ is $\alpha_{1}$, so
$a_{1}=(\omega,\alpha_{1}^{\vee})=1$ (see \S\ref{general}).
Moreover, $Z(L)^{0}/(Z(L)^{0}\cap H)$ is the one-parameter subgroup
of $T/T\cap H$ corresponding to $-3\omega_{2}^{\vee}$; indeed $(
-\omega_{2} , \omega_{2}^{\vee})<0$ and $ -\omega_{2} $ belongs to
the dual cone $\sigma^{\vee}$ associated to the closure of $T\cdot
x_{0}$ in $Z$. Since $Z(L)^{0}/(Z(L)^{0}\cap H)\equiv Z(L)^{0}\cdot
x_{0}$ is contained in $Z$, we have $1=(\omega,-3\omega_{2}^{\vee})=
-1-2a_{2}$, so
$\omega$ is $\omega_{1}-\omega_{2}$. We can study he simple variety
corresponding to the colored cone $(\sigma(\alpha_{2}^{\vee},
-\omega_{1}^{\vee}-\omega_{2}^{\vee}),$ $
\{D_{\alpha_{2}^{\vee}}\})$ in a similar way. Thus $Aut^{0}(X)$ is
isomorphic to $G$ by the Remark \ref{aut0}. $\square$

\begin{prop}\label{sl3/so3}
Suppose that $G/H$ is $SL_{3}/SO_{3}$ (and the type of $(G,\theta)$
is $AI$). Then $X$ is an hyperplane section of $\mathbb{LG}_{3}(6)$.
The full automorphism group $Aut(X)$ is generated by $SL_{3}$
($=Aut^{0}(X)$) and by $\theta$; moreover $Aut(X)$ has three orbits
in $X$.
\end{prop}

\begin{prop}\label{sl3}(see \cite{Bu}, Theorem 1.4)
Suppose that $G/H$ is $SL_{3}$. Then $X$ is an hyperplane section of
$\mathbb{G}_{3}(6)$. Moreover $Aut^{0}(X)$  is isomorphic to the
quotient of $SL_{3}\times SL_{3}$ by the intersection of the center
with the diagonal. The full automorphism group $Aut(X)$ is generated
by $Aut^{0}(X)$, by $\theta$ and by $(\varphi,\varphi)$ where
$\varphi$ is the automorphism of $SL_{3}$ corresponding to the
automorphism of the Dynkin diagram; in particular $Aut^{0}(X)$ has
index 4 in $Aut(X)$. Moreover $Aut(X)$ has three orbits in $X$.
\end{prop}

\begin{prop}\label{sl6/sp6}
Suppose that $G/H$ is $SL_{6}/Sp_{6}$ (and the involution's type is
$AII$). Then $X$ is an hyperplane section of $\mathbb{S}_{12}$. The
full automorphism group $Aut(X)$ is generated by $SL_{6}/\{\pm id\}$
($=Aut^{0}(X)$) and by $\theta$; thus $Aut^{0}(X)$ has index 2 in
$Aut(X)$. Moreover $Aut(X)$ has three orbits in $X$.
\end{prop}

\begin{prop}\label{E6/F4}
Suppose that $G/H$ is $E_{6}/F_{4}$ (and the involution's type is
$AII$). Then $X$ is an hyperplane section of $E_{7}/P_{7}\equiv
\mathbb{LG}(\mathbb{A}^{3},\mathbb{A}^{6})$. Moreover $Aut^{0}(X)$
is isomorphic to the simply-connected group of type $E_{6}$, while
the full automorphism group $Aut(X)$ is generated by $E_{6}$ and by
$\theta$; in particular  $Aut^{0}(X)$ has index 2 in $Aut(X)$.
Besides $Aut(X)$ has three orbits in $X$.
\end{prop}

Now, we study the automorphism group of $X$. By the Lemma
\ref{lemA2}, we know that $Lie(Aut(X))$ is $Lie(G)$. Notice that the
center of $G$ acts no-trivially over $X$. Moreover, if $G$ is
different from $SL_{3}\times SL_{3}$ and $SL_{6}$, then the center
of $G$ is a simple group, thus $Aut^{\circ}(X)$ is $G$. If $G$ is
$SL_{3}\times SL_{3}$, then the center is $C_{3}\times C_{3}$, where
$C_{n}$ is the group of $n$-th roots of the unit. Thus the
intersection of the center of $G$ with $G^{\theta}$ is the diagonal
of $C_{3}\times C_{3}$. If $G$ is $SL_{6}$ then the center is
$C_{6}$ and its intersection with $Sp_{6}(=G^{\theta})$ is $\{\pm
id\}$. Observe that $N_{G}(H)/H$ is simple because the fundamental
group of $R_{G,\theta}$ is $A_{2}$.

Now, we determine $K$ (see \S\ref{general}). It sufficient to
determine the subset $K'$ of $E$. Let $\varpi$ be the longest
element of $W$ and let $\varpi_{0}$ be the longest element of the
Weyl group of $R_{G}^{0}$ (here we consider the root system of $G$
with respect to $T$). Then $\theta':=\omega\omega_{0}\theta$ fixes
$T$ and $B$.  One can easily show that $\theta'$ exchange
$\omega_{1}$ with $\omega_{2}$ (see \cite{T2}, table 5.9 and
\cite{dCP1}, \S1.4), thus it exchanges $E_{1}$ with $E_{2}$ (recall
that the restriction of the valuation of $E_{i}$ to
$\mathbb{C}(X)^{(B)}/\mathbb{C}^{*}$ is $-\omega_{i}$). Hence, also
$\theta$ exchanges the previous two divisors; in particular $\theta$
exchanges the two closed orbits of $\widetilde{X}$. Notice that
$\theta$ belongs to $K'$ and induces the non-trivial automorphism of
the Dynkin diagram of $G$, with respect to $T'$ and $B'$ (see page
493 in \cite{LM} and \S26 in \cite{T2}). If $G$ is different from
$SL_{3}\times SL_{3}$, then $E/\overline{T}\,'$  contains exactly
two elements. Thus $K'/N_{T'}(H)$ coincides with
$E/\,\overline{T}\,'$. In particular, $Aut(G/H)$ is generated by
$Aut^{0}(G/H)$ and $\theta$. Moreover, $Aut(G/H)$ coincides with
$Aut(X)$ by Remark \ref{autX}.

Finally, suppose $G=SL_{3}\times SL_{3}$. Let $\dot{T}$ be a maximal
torus of $SL_{3}$, let $\dot{B}$ be a Borel subgroup of $SL_{3}$ and
let $\varphi$ be the equivariant automorphism of $SL_{3}$ associated
to the non-trivial automorphism of the Dynkin diagram (with respect
to $\dot{T}$ and $\dot{B}$). We can set $T=T'=\dot{T}\times
\dot{T}$, $B=\dot{B}\times \dot{B}^{\,\text{-}}$ and
$B'=\dot{B}\times \dot{B}$.  One can easily see that $K'$ is
generated by $N_{T'}(H)$, $\theta$ and $(\varphi,\varphi)$. (Notice
that $E/\overline{T}\,'$ has eight elements: $id$, $\theta$,
$(\varphi,id)$, $(id,\varphi)$, $(\varphi,\varphi)$,
$\theta\circ(\varphi,\varphi)$, $\theta\circ(\varphi,id)$ and
$\theta\circ(id,\varphi)$). Notice that $(\varphi,\varphi)$
stabilizes both $B$ and $B^{\,\text{-}}$. We have
$\omega_{i}=\varpi_{i}^{1}-\varpi_{i}^{2}$ where
$\{\varpi_{1}^{j},\varpi_{2}^{j}\}$ are the fundamental weights of
the $j$-th copy of $SL_{3}$ in $G$ (with respect to $\dot{T}$ and
$\dot{B}$). Thus $(\varphi,\varphi)$ exchange $\omega_{1}$ with
$\omega_{2}$, so $(\varphi,\varphi)$ exchange $E_{1}$ with $E_{2}$.
Furthermore, $ (\varphi,\varphi)$ induces a non-trivial automorphism
of the Dynkin diagram of $G$ (with respect to $T'$ and $B'$).
Therefore, $Aut(X)$ coincides with $Aut(G/H)$ by Remark \ref{autX}.

\section{Varieties of rank at least three}\label{rank3}

{\em In the following we always suppose that the rank of $G/H$ is at
least two.} In this section we consider the symmetric varieties
which belongs to an infinite family; in particular, we consider all
the symmetric varieties of rank at least three. We consider also
completions of  the following varieties of rank two: 1) the
symmetric variety $PSL_{3}$; 2) the symmetric variety
$PSL_{3}/PSO_{3}$ of type $AI$; 3) the symmetric variety
$PSL_{4}/PSp_{4}$ of type $AII$; 4) the symmetric variety $SO_{5}$;
5) the symmetric variety $Spin_{5}\equiv Sp_{4}$; 6) the symmetric
varieties $SO_{n}/N_{SO_{n}}(SO_{2}\times SO_{m-2})$ of type $BI$
and $DI$; 7) the symmetric varieties
$Sp_{2n}/N_{Sp_{2n}}(Sp_{2}\times Sp_{2n-2})$ of type $CII$; 8) the
symmetric variety $SO_{8}/N_{SO_{8}}(GL_{4})$ of type $DIII$.

Given a linear endomorphism $\varphi$ of a vector space $V$, let
$graph(\varphi)$ be the subspace $\{(v,\varphi(v)): v\in V\}$ of
$V\oplus V$. Observe that $graph(\varphi)$ has the same dimension of
$V$. If we have fixed  a (skew-)symmetric bilinear form on $V$, then
we define a (skew-)symmetric bilinear $q'$ form on $V\oplus V$ such
that $q'(v,w)=q(v)-q(w)$ for each $(v,w)\in V\oplus V$.

\subsection{Restricted root system of type $A_{l}$}\label{cap-Al}

Now, we consider the symmetric varieties such that: 1)
$H=N_{G}(G^{\theta})$; 2) the  restricted root system has type
$A_{l}$ (with $l\geq2$).

\begin{prop}\label{Al}
Let $X$ be a smooth projective symmetric variety with Picard number
one and rank at least two. Suppose that $H=N_{G}(G^{\theta})$ and
that the restricted root system has type $A_{l}$, then $X$ is
isomorphic to the projectivization of an irreducible
$G$-representation.
\end{prop}

{\em Proof.} It is sufficient to show an irreducible spherical
representation with dimension equal to $\textrm{dim}\ G/H+1$. We
have already considered the case of $E_{6}/N_{E_{6}}(F_{4})$ in
Proposition \ref{EIV}.

1) If $G/H$ is isomorphic to $PGL_{l+1} $, then $X$ is isomorphic to
$\mathbb{P}(M_{l+1})$ as $(PGL_{l+1}\times PGL_{l+1})$-variety (here
$M_{l+1} $ is the space of matrices of order $l+1$). Indeed
$PGL_{l+1} \times PGL_{l+1} $ acts on $ M_{l+1} $ and the stabilizer
of $\mathbb{C}\textit{1}  $ is the diagonal, namely $G^{\theta}$.
Moreover $\mathbb{P}(M_{l+1} )$ has dimension equal to $l^{2}+l$.

2) If $G/H$ is isomorphic to $SL_{l+1} /N_{SL_{l+1} }(SO_{l+1} )$
and $(G,\theta)$ has type $AI$,  then $X$ is isomorphic to
$\mathbb{P}(Sym^{2}(\mathbb{C}^{l+1})) $
(see \cite{CM}, Theorem 2.3 or \cite{T2}, Proposition 26.4).

3) If $G/H$ is isomorphic to $SL_{l+1} /N_{SL_{l+1} }(Sp_{l+1} )$
and $(G,\theta)$ has type $AII$,  then $X$ is isomorphic to
$\mathbb{P}(\bigwedge^{2l-2}(\mathbb{C}^{2l})) $.
$\square$

\subsection{Restricted root system of type $B_{l}$ or
$D_{l}$}\label{cap Bl+Dl}

Now we consider the symmetric varieties whose restricted root system
has type $D_{l}$ or $B_{l}$. In the second case we suppose also
$H=N_{G}(G^{\theta})$. Explicitly, we consider the following cases:
1) the symmetric varieties $SO_{n}$; 2) the symmetric varieties
$SO_{n}/(SO_{l}\times SO_{n-l})$ of type $BI$ and $DI$. We consider
also the symmetric variety $SO_{8}/GL_{4}$ of type $DIII$ because it
is isomorphic to $SO_{8}/(SO_{2}\times SO_{6})$.

\begin{prop}\label{Bl+Dl}
The smooth projective embedding of $SO_{n}$ with Picard number one
is isomorphic to $\mathbb{S}_{n}$.
\end{prop}

{\em Proof.} We have an inclusion of $SO_{n}$ into
$\mathbb{I}\mathbb{G}_{n}(2n)$ given by the map $g\rightarrow
graph(\varphi)$. It is easy to show that this map is compatible with
the action of $SO_{n}\times SO_{n}$ on $SO_{n}$ and
$\mathbb{I}\mathbb{G}_{n}(2n)$. Thus it is sufficient to observe
that $SO_{n}$ and $\mathbb{I}\mathbb{G}_{n}(2n)$ have the same
dimension. $\square$

\begin{prop}\label{BI+DI}
The smooth projective embedding of $SO_{n}/(SO_{l}\times SO_{n-l})$
with Picard number one is isomorphic to $\mathbb{G}_{l}(n)$.
\end{prop}

{\em Proof.} Let $W$ be a $l$-dimensional subspace of
$\mathbb{C}^{n}$ over which the quadratic form is nondegenerate. We
denote by $W^{\bot}$ the orthogonal complement of $W$. The
stabilizer of $W$ contains $SO(W)\times SO(W^{\bot})$. Thus
$\mathbb{G}_{l}(n)$ is a smooth completion with Picard number one of
the symmetric variety $SO_{n}/Stab_{SO_{n}}(W)$ (see \cite{dCP1},
Lemma 1.7).
This fact implies, by Theorem \ref{classif}, that $H$ is
$SO(W)\times SO(W^{\bot})$ (notice that $G/H$ cannot be
$Spin_{n}/(Spin_{2}\times Spin_{n-2})$ because this last variety is
Hermitian and its restricted root system has type $B_{2}$).
%
%
%
%
%
%
%
%
$\square$

\begin{prop}\label{DIII}
The smooth projective embedding of $SO_{8}/GL_{4}$ with Picard
number one is isomorphic to $\mathbb{G}_{2}(8)$.
\end{prop}

{\em Proof.} Observe that $SO_{8}/GL_{4}$ is isomorphic to
$SO_{8}/S(O_{2}\times O_{6})$. Indeed the involution of $SO_{8}$ of
type $DIII$ (and rank two) is conjugated by an equivariant
automorphism of $SO_{8}$ to  the involution of $SO_{8}$ of type $DI$
and (and rank two). Indeed, these varieties have the same  Satake
diagrams up to an automorphism of the Dynkin diagram of $SO_{8}$.
Now, it is sufficient to observe that a homogeneous symmetric
variety has at most one smooth completion with Picard number one.
$\square$

\subsection{Restricted root system with type $BC_{l}$ or $
C_{l}$}\label{cap BCl+Cl}

Now, we consider the symmetric varieties whose restricted root
system has type $BC_{l}$ or $ C_{l}$. We have the following two
cases: 1) the symmetric varieties $Sp_{2n}/(Sp_{2l}\times
Sp_{2n-2l})$ of type $CII$; 2) the symmetric varieties $Sp_{2l}$. We
consider also $Spin_{5}$, because it is isomorphic to $Sp_{4}$.

\begin{prop}\label{CII}
The smooth projective embedding of $Sp_{2n}/Sp_{2l}\times
Sp_{2n-2l}$ with Picard number one is isomorphic to
$\mathbb{G}_{2l}(2n)$.
\end{prop}

{\em Proof.} Let $W$ be a $2l$-dimensional subspace of
$\mathbb{C}^{2n}$ over which the standard bilinear skew-symmetric
form is nondegenerate and let $W^{\bot}$ be its orthogonal
complement. The $Sp_{2n}$-orbit of $W$ in $\mathbb{G}_{2l}(2n)$ is
isomorphic to $Sp_{2n}/Sp(W)\times Sp(W^{\bot})$ (see \cite{dCP1},
Lemma 1.7), thus is dense in $\mathbb{G}_{2l}(2n)$. $\square$

\begin{prop}\label{Cl}
The smooth projective embedding of $Sp(2l)$ with Picard number one
is isomorphic to $\mathbb{I}\mathbb{G}_{2l}(4l)$.
\end{prop}

{\em Proof.} We have an inclusion of $Sp(2l)$ into
$\mathbb{I}\mathbb{G}_{2l}(4l)$ given by the map $g\rightarrow
graph(\varphi)$. It is easy to show that this map is compatible with
the action of $Sp(2l)\times Sp(2l)$ over $Sp(2l)$ and
$\mathbb{I}\mathbb{G}_{2l}(4l)$. Furthermore, $Sp(2l)$ and
$\mathbb{I}\mathbb{G}_{2l}(4l)$ have the same dimension. $\square$

\begin{prop}\label{C2}
The smooth projective embedding of $Spin_{5}$ with Picard number one
is isomorphic to $\mathbb{L}\mathbb{G}_{4}(8)$.
\end{prop}

{\em Proof.} Observe that $Spin_{5}\times Spin_{5}$ is isomorphic to
$Sp_{4}\times Sp_{4}$. $\square$

\section*{Acknowledgments}

We would like to thank M. Brion  for the continued support and for
many very helpful suggestions. I would also thank L. Manivel for
many enlightening suggestion. The author has been partially
supported by C.N.R.S.'S grant in collaboration with the Liegrits.

\end{document}